%% file: IV-5-CR-equivalences.tex
\documentclass[12pt,twoside,leqno,openany]{amsart}

\usepackage{amssymb,amsbsy,amsmath,amsfonts,amssymb,amscd,times,
graphics,color,xypic,footmisc,fancyhdr,multicol,fancybox,
graphicx,mathrsfs,rotating,ifthen,wasysym}
\usepackage[all]{xy}

\usepackage[T1]{fontenc}
\input macros.tex

\let\mathcal\mathscr

\begin{document}

$\:$

\bigskip\bigskip

\begin{center}

{\Large\bf Equivalences of $5$-dimensional CR-manifolds IV:}

\medskip

{\Large\bf Six ambiguity matrix groups}

\medskip

{\Large\bf (Initial $G$-structures)}

\end{center}

\medskip

\begin{center}
Jo\"el {\sc Merker}
\end{center}

\bigskip

\begin{center}
\begin{minipage}[t]{10.25cm}
\baselineskip =0.32cm 
{\scriptsize
{\bf Abstract.}
Class $\text{\sf I}$ CR manifolds have initial $G$-structure a certain
4-dimensional subgroup of ${\sf GL}_3(\C)$.  Class $\text{\sf II}$ CR
manifolds have initial $G$-structure a certain 10-dimensional subgroup
of ${\sf GL}_4(\C)$. Class $\text{\sf III}_{\text{\sf 1}}$ CR
manifolds have initial $G$-structure a certain 10-dimensional subgroup
of ${\sf GL}_5(\C)$. Class $\text{\sf III}_{\text{\sf 2}}$ CR
manifolds have initial $G$-structure a certain 18-dimensional subgroup
of ${\sf GL}_5(\C)$. Class $\text{\sf IV}_{\text{\sf 1}}$ CR
manifolds have initial $G$-structure a certain 13-dimensional subgroup
of ${\sf GL}_5(\C)$. Class $\text{\sf IV}_{\text{\sf 2}}$ CR
manifolds have initial $G$-structure a certain 10-dimensional subgroup
of ${\sf GL}_5(\C)$.
}
\end{minipage}
\end{center}

\medskip

\begin{center}
\begin{minipage}[t]{11.75cm}
\baselineskip =0.35cm {\scriptsize

\centerline{\bf Table of contents}

\medskip

{\bf \ref{general-class-I}.~General class $\text{\sf I}$: 
$4$-dimensional real subgroup of ${\sf GL}_3(\C)$
\dotfill~\pageref{general-class-I}.}

{\bf \ref{multiplicator-a}.~Scholium: The multiplicator function $a$
\dotfill~\pageref{multiplicator-a}.}

{\bf \ref{general-class-II}.~General class $\text{\sf II}$
\dotfill~\pageref{general-class-II}.}

{\bf \ref{general-class-III-1}.~General class 
$\text{\sf III}_{\text{\sf 1}}$
\dotfill~\pageref{general-class-III-1}.}

{\bf \ref{general-class-III-2}.~General class 
$\text{\sf III}_{\text{\sf 2}}$
\dotfill~\pageref{general-class-III-2}.}

{\bf \ref{general-class-IV-1}.~General class 
$\text{\sf IV}_{\text{\sf 1}}$
\dotfill~\pageref{general-class-IV-1}.}

{\bf \ref{general-class-IV-2}.~General class 
$\text{\sf IV}_{\text{\sf 2}}$
\dotfill~\pageref{general-class-IV-2}.}

}\end{minipage}
\end{center}


\bigskip

\section{\sf General class $\text{\sf I}$: 
$4$-dimensional real subgroup of ${\sf GL}_3(\C)$}
\label{general-class-I}
\HEAD{\ref{general-class-I}.~General class $\text{\sf I}$: 
$4$-dimensional real subgroup of ${\sf GL}_3(\C)$}{
Jo\"el {\sc Merker}, D\'epartement de Math\'ematiques d'Orsay}

\medskip

Equip $\C^2$ with coordinates:
\[
(z,w)
\,\in\,
\C^2.
\]

Let a connected real hypersurface:
\[
M^3
\subset
\C^2
\]
be of smoothness:
\[
\mathcal{C}^\kappa\ \
(\kappa\geqslant 3),
\ \ \ \ \ \ \ \ \ \ \ \ \ \ \ \ \ \ \ \ \
\text{\rm or}
\ \ \ \ \ \ \ \ \ \ \ \ \ \ \ \ \ \ \ \ \
\mathcal{C}^\infty,
\ \ \ \ \ \ \ \ \ \ \ \ \ \ \ \ \ \ \ \ \
\text{\rm or}
\ \ \ \ \ \ \ \ \ \ \ \ \ \ \ \ \ \ \ \ \
\mathcal{C}^\omega.
\]

Pick a point:
\[
p\in M,
\]
and take a (small) open neighborhood:
\[
p
\in
{\sf U}_p
\subset
\C^2.
\]

By definition (\cite{ Merker-Pocchiola-Sabzevari-5-CR-II,
Merker-5-CR-III}):
\[
\Big(
M^3
\,\subset\,
\C^2
\Big)
\,\,\in\,\,
\text{\sf General Class $\text{\sf I}$},
\]
belongs to the first class if:
\[
\C\otimes_\R TM
=
T^{1,0}M
+
T^{0,1}M
+
\big[
T^{1,0}M,\,T^{0,1}M\big].
\]

This means that for any local vector field generator:
\[
\mathcal{L}
\,=\,
\text{\rm section of}\,\,
T^{1,0}\big(M\cap{\sf U}_p\big),
\]
one has at every point $q \in M\cap {\sf U}_p$:
\[
{\bf 3}
=
\rank_\C
\Big(
\mathcal{L}\big\vert_q,\,\,
\overline{\mathcal{L}}\big\vert_q,\,\,
\big[\mathcal{L},\overline{\mathcal{L}}\big]
\Big\vert_q
\Big).
\]

\medskip

Next, take any (local) biholomorphism:
\[
h\colon\ \ \
{\sf U}_p
\overset{\sim}{\,\longrightarrow\,}
{\sf U}_{p'}'
=
h({\sf U}_p)
\ \ \ \ \ \ \ \ \ \ \ \ \
{\scriptstyle{(p'\,=\,\,h(p))}},
\]
which, when ${\sf U}_p$ is small enough, certainly 
transfers $M \cap {\sf U}_p$ to a certain hypersurface:
\[
{M'}^3
\,:=\,
h\big(M\cap{\sf U}_p\big)
\,\subset\,
{\C'}^2,
\]
having the same smoothness (mental exercise):
\[
\mathcal{C}^\kappa\ \
(\kappa\geqslant 3),
\ \ \ \ \ \ \ \ \ \ \ \ \ \ \ \ \ \ \ \ \
\text{\rm or}
\ \ \ \ \ \ \ \ \ \ \ \ \ \ \ \ \ \ \ \ \
\mathcal{C}^\infty,
\ \ \ \ \ \ \ \ \ \ \ \ \ \ \ \ \ \ \ \ \
\text{\rm or}
\ \ \ \ \ \ \ \ \ \ \ \ \ \ \ \ \ \ \ \ \
\mathcal{C}^\omega.
\]

\medskip

Take also any local vector field generator:
\[
\mathcal{L}'
\,=\,
\text{\rm section of}\,\,
T^{1,0}M'.
\]
Then necessarily at every point $q' \in M'$, one 
also has (\cite{ Merker-Pocchiola-Sabzevari-5-CR-II}, Section~4):
\[
{\bf 3}
=
\rank_\C
\Big(
\mathcal{L}'\big\vert_{q'},\,\,
\overline{\mathcal{L}}'\big\vert_{q'},\,\,
\big[\mathcal{L}',\overline{\mathcal{L}}'\big]
\Big\vert_{q'}
\Big).
\]

Write the components of $h$ and the target coordinates as:
\[
\aligned
(z,w)
&
\,\longmapsto\,
\big(z'(z,w),\,w'(z,w)\big)
\\
&\ \ \ \
=:
(z',w').
\endaligned
\]

To lower the number of primes in equations, it proves
to be better to work instead with the {\em inverse}
of $h$, denoted here not with an exponent ${}^{-1}$,
but with a {\em prime:}
\[
h'\colon \ \ \ \ \ \ \ \ \ \ \
\aligned
{\sf U}_{p'}'
&
\overset{\sim}{\,\longrightarrow\,}
{\sf U}_p
\ \ \ \ \ \ \ \ \ \ \ \ \ \ \ \ \ \ \ \ \ \ \ \ \ \ \ \ \ \
{\scriptstyle{(p\,=\,\,h'(p'))}},
\\
(z',w')
&
\,\longmapsto\,
\big(z(z',w'),\,w(z',w')\big)
=
(z,w).
\endaligned
\]

Then from~\cite{ Merker-Pocchiola-Sabzevari-5-CR-II}, there
exists a nowhere vanishing function:
\[
a\colon\ \ \
M\cap{\sf U}_p
\,\longrightarrow\,
\C\backslash\{0\},
\]
such that:
\[
h_*'
\big(\mathcal{L}'\big)
=
a\,\mathcal{L}.
\]

Of course simultaneously (\cite{ Merker-Pocchiola-Sabzevari-5-CR-II}):
\[
h_*'
\big(\overline{\mathcal{L}}'\big)
=
\overline{a}\,\overline{\mathcal{L}}.
\]

Now set:
\[
\aligned
\mathcal{T}
&
:=
\isqrt\,
\big[\mathcal{L},\,\overline{\mathcal{L}}\big],
\\
\mathcal{T}'
&
:=
\isqrt\,
\big[\mathcal{L}',\,\overline{\mathcal{L}}'\big].
\endaligned
\]
The advantage of the $\isqrt$-factor is {\em reality:}
\[
\aligned
\overline{\mathcal{T}}
&
=
\mathcal{T},
\\
\overline{\mathcal{T}}'
&
=
\mathcal{T}',
\endaligned
\]
so that one has two (local) frames:
\[
\aligned
&
\big\{\mathcal{L},\,\overline{\mathcal{L}},\,
\mathcal{T}\big\}
\ \ \ \ \ \ \ \,
\text{\rm for}\ \
\C\otimes_\R T(M\cap{\sf U}_p),
\\
&
\big\{\mathcal{L}',\,\overline{\mathcal{L}}',\,
\mathcal{T}'\big\}
\ \ \ \ \
\text{\rm for}\ \
\C\otimes_\R TM'.
\endaligned
\]

The Lie bracket $\mathcal{ T}'$ transfers
back to $M \cap {\sf U}_p$ through $h_*'$ as:
\[
\aligned
h_*'\big(\mathcal{T}'\big)
&
=
h_*'\big(\isqrt\,\big[\mathcal{L}',\,\overline{\mathcal{L}}'\big]\big)
\\
&
=
\isqrt\,
\big[
h_*'\big(\mathcal{L}'\big),\,
h_*'\big(\overline{\mathcal{L}}'\big)
\big]
\\
&
=
\isqrt\,
\big[a\,\mathcal{L},\,\,
\overline{a}\,\overline{\mathcal{L}}\big]
\\
&
=
a\overline{a}\,
\underbrace{\isqrt\,\,\big[\mathcal{L},\,\overline{\mathcal{L}}\big]}_{
=\,\mathcal{T}}
+
\underbrace{
\isqrt\,a\,\mathcal{L}\big(\overline{a}\big)}_{
=:\,\overline{b}}
\cdot
\overline{\mathcal{L}}
\underbrace{-
\isqrt\,\overline{a}\,\overline{\mathcal{L}}(a)}_{
=:\,b}
\cdot
\mathcal{L},
\endaligned
\]
so that, if one decides to set:
\[
b
:=
-\,\isqrt\,\overline{a}\,\overline{\mathcal{L}}(a),
\]
forgetting how this coefficient-function
is related to $a$, one obtains:
\[
h_*'\big(\mathcal{T}'\big)
=
a\overline{a}\,\mathcal{T}
+
\overline{b}\,\overline{\mathcal{L}}
+
b\,\mathcal{L}.
\]

\medskip\noindent{\bf Summary.}
{\em Through any local biholomorphic equivalences
between hypersurfaces of $\C^2$ belonging
to the General Class $\text{\sf I}$:}
\[
M^3
\overset{\sim}{\,\longrightarrow\,}
{M'}^3,
\]
{\em for any two choices of local vector
field generators:}
\[
\aligned
&
\mathcal{L}
\ \ \ \ \
\text{\rm for}\ \
T^{1,0}M,
\\
&
\mathcal{L}'
\ \ \ \ \
\text{\rm for}\ \
T^{1,0}M',
\endaligned
\]
{\em the transfer of frame obeys the rule:}
\[
\left(\!\!
\begin{array}{c}
\mathcal{L}'
\\
\overline{\mathcal{L}}'
\\
\mathcal{T}'
\end{array}
\!\!\right)
=
\left(\!
\begin{array}{ccc}
a & 0 & 0
\\
0 & \overline{a} & 0
\\
b & \overline{b} & a\overline{a}
\end{array}
\!\right)
\left(\!\!
\begin{array}{c}
\mathcal{L}
\\
\overline{\mathcal{L}}
\\
\mathcal{T}
\end{array}
\!\!\right),
\]
{\em for some two local functions:}
\[
a\colon\ \ \
M
\,\longrightarrow\,
\C\backslash\{0\},
\ \ \ \ \ \ \ \ \ \ \ \ \ \ \ \ \ \ \ \ \ \
b\colon\ \ \
M
\,\longrightarrow\,
\C.
\]

\medskip

This means that the {\em ambiguity} in the choice of a local frame:
\[
\big\{
\mathcal{L},\,
\overline{\mathcal{L}},\,
\mathcal{T}
\big\}
\ \ \ \ \
\text{\rm for}\ \
\C\otimes_\R TM
\]
which comes\,\,---\,\,{\em naturally from the
point of view of CR geometry}\,\,---\,\,from a choice of a local frame:
\[
\mathcal{L}
\ \ \ \ \
\text{\rm for}\ \
T^{1,0}M
\]
is represented by general changes of frames whose
matrices are of the form: 
\[
\left(\!
\begin{array}{ccc}
a & 0 & 0
\\
0 & \overline{a} & 0
\\
b & \overline{b} & a\overline{a}
\end{array}
\!\right),
\]
the entries being coefficient-functions depending on existing
equivalences $M \overset{ \sim}{ \longrightarrow} M'$,
and\big/or on the choice of local coordinates.

\medskip

Within \'Elie Cartan's theory of equivalence between
differential-geometric structures, this means that the {\em initial
$G$-structure} for the biholomorphic equivalence problem between
hypersurface $M^3 \subset \C^2$ belonging to the General Class
$\text{\sf I}$ is a reduction of the full linear group 
${\sf GL}_3 ( \C)$ to
the mentioned subgroup in which the\,\,---\,\,possibly
unknown\,\,---\,\,functions $a$, $b$ are {\em replaced} by independent
complex variables.

\medskip

Indeed, one has a very elementary:

\medskip\noindent{\bf Lemma.}
{\em The set of matrices:}
\[
{\sf G}_{\text{\sf I}}^{\sf initial}
\,:=\,
\left\{
\left(\!\!
\begin{array}{ccc}
{\sf a} & 0 & 0
\\
0 & \overline{\sf a} & 0
\\
{\sf b} & \overline{\sf b} & {\sf a}\overline{\sf a}
\end{array}
\!\!\right)
\,\in\,
\mathcal{M}_{3\times 3}(\C)\,\colon\,\,
{\sf a}\,\in\,\C\backslash\{0\},\,\,
{\sf b}\,\in\,\C
\right\}
\]
{\em is a (closed) $4$-dimensional
real matrix subgroup of the full:}
\[
{\sf GL}_3(\C)
\,:=\,
\left\{
\pi
=
\left(\!\!
\begin{array}{ccc}
\pi_{1,1} & \pi_{1,2} & \pi_{1,3}
\\
\pi_{2,1} & \pi_{2,2} & \pi_{3,3}
\\
\pi_{3,1} & \pi_{3,2} & \pi_{3,3}
\end{array}
\!\!\right)
\,\in\,
\mathcal{M}_{3\times 3}(\C)\,\colon\,\,
0
\neq
\det\,\pi
\right\}.
\]

\proof
Well, closedness under multiplication (composition):
\[
\aligned
\left(\!\!
\begin{array}{ccc}
{\sf a}_1 & 0 & 0
\\
0 & \overline{\sf a}_1 & 0
\\
{\sf b}_1 & \overline{\sf b}_1 &
{\sf a}_1\overline{\sf a}_1
\end{array}
\!\!\right)
\cdot
\left(\!\!
\begin{array}{ccc}
{\sf a}_2 & 0 & 0
\\
0 & \overline{\sf a}_2 & 0
\\
{\sf b}_2 & \overline{\sf b}_2 &
{\sf a}_2\overline{\sf a}_2
\end{array}
\!\!\right)
&
=
\left(\!\!
\begin{array}{ccc}
{\sf a}_1{\sf a}_2 & 0 & 0
\\
0 & \overline{\sf a}_1\overline{\sf a}_2 & 0
\\
{\sf b}_1{\sf a}_2
\!+\!{\sf a}_1\overline{\sf a}_1{\sf b}_2 & 
\overline{\sf b}_1\overline{\sf a}_2
\!+\!{\sf a}_1\overline{\sf a}_1\overline{\sf b}_2 &
{\sf a}_1\overline{\sf a}_1{\sf a}_2\overline{\sf a}_2
\end{array}
\!\!\right)
\\
&
=:
\left(\!\!
\begin{array}{ccc}
{\sf a}_3 & 0 & 0
\\
0 & \overline{\sf a}_3 & 0
\\
{\sf b}_3 & \overline{\sf b}_3 &
{\sf a}_3\overline{\sf a}_3
\end{array}
\!\!\right)
\endaligned
\]
is visibly clear, after setting:
\[
\aligned
{\sf a}_3
&
:=
{\sf a}_1{\sf a}_2
\ \ \ \ \ \ \ \ \ \ \ \ \ \ \ \ \ \ \ \ \ \ \ \ \
{\scriptstyle{(\neq\,0,\,\,
\text{\rm again})}}
\\
{\sf b}_3
&
:=
{\sf b}_1{\sf a}_2
+
{\sf a}_1\overline{\sf a}_1{\sf b}_2.
\endaligned
\]

Quite similarly, the inverse:
\[
\aligned
\left(\!\!
\begin{array}{ccc}
{\sf a} & 0 & 0
\\
0 & \overline{\sf a} & 0
\\
{\sf b} & \overline{\sf b} &
{\sf a}\overline{\sf a}
\end{array}
\!\!\right)^{-1}
&
=
\left(\!\!
\begin{array}{ccc}
\frac{1}{\sf a} & 0 & 0
\\
0 & \frac{1}{\overline{\sf a}} & 0
\\
\frac{-{\sf b}}{{\sf a}{\sf a}\overline{\sf a}} & 
\frac{-\overline{\sf b}}{{\sf a}\overline{\sf a}\overline{\sf a}} &
\frac{1}{{\sf a}\overline{\sf a}}
\end{array}
\!\!\right)
\\
&
=:
\left(\!\!
\begin{array}{ccc}
{\sf a}^\sim & 0 & 0
\\
0 & \overline{\sf a}^\sim & 0
\\
{\sf b}^\sim & \overline{\sf b}^\sim &
{\sf a}^\sim\overline{\sf a}^\sim
\end{array}
\!\!\right),
\endaligned
\]
also belongs to the subgroup, with:
\[
\aligned
{\sf a}^\sim
&
:=
\frac{1}{\sf a}
\ \ \ \ \ \ \ \ \ \ \ \ \ \ \ \ \ \ \ \ \ \ \ \ \
{\scriptstyle{(\neq\,0)}},
\\
{\sf b}^\sim
&
:=
-\,\frac{\sf b}{{\sf a}{\sf a}\overline{\sf a}},
\endaligned
\]
which concludes.
\endproof

\noindent{\bf Proposition.}
{\em On a $3$-dimensional hypersurface:}
\[
\Big(
M^3
\,\subset\,
\C^2
\Big)
\,\,\in\,\,
\text{\sf General Class $\text{\sf I}$},
\]
{\em having biholomorphically invariant $(1, 0)$ CR bundle:}
\[
T^{1,0}M
\,\subset\,
\C\otimes_\R TM,
\]
{\em for any choice of local vector field generator:}
\[
\mathcal{L}
\ \ \ \ \
\text{\rm for}\ \
T^{1,0}M,
\]
{\em the associated frame:}
\[
\big\{
\mathcal{L},\,\overline{\mathcal{L}},\,
\isqrt\,\big[\mathcal{L},\overline{\mathcal{L}}\big]
\big\}
\,=:\,
\big\{
\mathcal{L},\,\overline{\mathcal{L}},\,\mathcal{T}
\big\}
\]
{\em for the full complexified tangent bundle:}
\[
\C\otimes_\R TM
\]
{\em performs a reduction of the full ${\sf GL}_3 ( \C)$-structure
of $\C \otimes_\R TM$ to the $4$-dimensional subgroup:}
\[
{\sf G}_{\text{\sf I}}^{\sf initial}
\,:=\,
\left\{
\left(\!\!
\begin{array}{ccc}
{\sf a} & 0 & 0
\\
0 & \overline{\sf a} & 0
\\
{\sf b} & \overline{\sf b} & {\sf a}\overline{\sf a}
\end{array}
\!\!\right)
\,\in\,
\mathcal{M}_{3\times 3}(\C)\,\colon\,\,
{\sf a}\,\in\,\C\backslash\{0\},\,\,
{\sf b}\,\in\,\C
\right\},
\]
{\em which is the initial ambiguity group when
launching the Cartan equivalence method.\qed}


\bigskip

\section{\sf Scholium: The multiplicator function $a$}
\label{multiplicator-a}
\HEAD{\ref{multiplicator-a}.~Scholium: The multiplicator function $a$}{
Jo\"el {\sc Merker}, D\'epartement de Math\'ematiques d'Orsay}

\medskip

Although Cartan's method equivalence classicaly {\em decides} not to
enter the computations of coefficient-functions
like $a$ and $b$ in the preceding section,
one may wonder what is the expression of $a$, at least. 

Through a (local) biholomorphism:
\[
h'
\colon\ \ \
(z',w')
\,\longmapsto\,
\big(z(z',w'),\,w(z',w')\big),
\]
basic vector fields transfer as
(without writing $h_*'$):
\[
\aligned
\frac{\partial}{\partial z'}
&
=
z_{z'}\,\frac{\partial}{\partial z}
+
w_{z'}\,\frac{\partial}{\partial w},
\\
\frac{\partial}{\partial w'}
&
=
z_{w'}\,\frac{\partial}{\partial z}
+
w_{w'}\,\frac{\partial}{\partial w}.
\endaligned
\]

If $M$ and $M'$ have local graphing equations:
\[
\aligned
v
&
=
\varphi(x,y,u),
\\
v'
&
=
\varphi'(x',y',u'),
\endaligned
\]
two generators for $T^{1, 0}M$ and $T^{1,0}M'$ are
(\cite{Merker-Pocchiola-Sabzevari-5-CR-II}):
\[
\aligned
\mathcal{L}
&
=
\frac{\partial}{\partial z}
+
2A\,
\frac{\partial}{\partial w},
\\
\mathcal{L}'
&
=
\frac{\partial}{\partial z'}
+
2A'\,
\frac{\partial}{\partial w'},
\endaligned
\]
where:
\[
\aligned
A
&
:=
\frac{\varphi_z}{\isqrt+\varphi_u},
\\
A'
&
:=
\frac{\varphi_{z'}'}{\isqrt+\varphi_{u'}'}.
\endaligned
\]

Consequently:
\[
\aligned
\mathcal{L}'
&
=
\frac{\partial}{\partial z'}
+
2A'\,\frac{\partial}{\partial w'}
\\
&
=
\big(z_{z'}+2A'\,z_{w'}\big)\,
\frac{\partial}{\partial z}
+
\big(
w_{z'}+2A'\,w_{w'}\big)\,
\frac{\partial}{\partial w}
\\
&
=
\big(
z_{z'}+2A'\,z_{w'}
\big)\,
\bigg[
\underbrace{\frac{\partial}{\partial z}
+
\frac{w_{z'}+2A'\,w_{w'}}{z_{z'}+2A'\,z_{w'}}\,
\frac{\partial}{\partial w}}_{
=\,\mathcal{L}}
\bigg],
\endaligned
\]
so that necessarily:
\[
\aligned
a
&
=
z_{z'}
+
2A'\,z_{w'},
\\
A
&
=
\frac{w_{z'}+2A'\,w_{w'}}{
z_{z'}+2A'\,z_{w'}},
\endaligned
\]
or more precisely writing all arguments:
\[
\aligned
a(z,w)
&
=
\big(z_{z'}+2A'\,z_{w'}\big)
\big(z'(z,w),\,w'(z,w)\big),
\\
A(z,w)
&
=
\frac{w_{z'}+2A'\,w_{w'}}{
z_{z'}+2A'\,z_{w'}}
\big(z'(z,w),\,w'(z,w)\big),
\endaligned
\]
so that:
\[
\aligned
h_*'\big(\mathcal{L}'\big)
&
=
h_*'\bigg(
\frac{\partial}{\partial z'}
+
A'\,\frac{\partial}{\partial w'}
\bigg)
\\
&
=
a\,
\bigg(
\frac{\partial}{\partial z}
+
A\,\frac{\partial}{\partial w}
\bigg)
\\
&
=
a\,\mathcal{L}.
\endaligned
\]

Here, because $h'$ is (local) biholomorphism,
its Jacobian matrix $h_*'$ is invertible,
hence on restriction to any
complex $1$-dimensional line of the form:
\[
T_{q'}^{1,0}M'
=
\C\cdot\mathcal{L}'\big\vert_{q'}
\ \ \ \ \ \ \ \ \ \ \ \ \
{\scriptstyle{(q'\,\in\,M')}},
\] 
it is of rank $1$ (also invertible):
\[
h_*'\colon\ \ \
T_{q'}^{1,0}M'
\overset{\sim}{\,\longrightarrow\,}
T_{h'(q')}^{1,0}M,
\]
from which it follows that the function $a$, 
appearing also in denominator place in the
above expression of $A$, vanishes {\em nowhere}.


\bigskip

\section{\sf General class $\text{\sf II}$}
\label{general-class-II}
\HEAD{\ref{general-class-II}.~General class $\text{\sf II}$}{
Jo\"el {\sc Merker}, D\'epartement de Math\'ematiques d'Orsay}

\medskip

Equip $\C^3$ with coordinates:
\[
(z,w_1,w_2)
\,\in\,
\C^3.
\]

Let a connected CR-generic submanifold:
\[
M^4
\subset
\C^3
\]
be of smoothness:
\[
\mathcal{C}^\kappa\ \
(\kappa\geqslant 3),
\ \ \ \ \ \ \ \ \ \ \ \ \ \ \ \ \ \ \ \ \
\text{\rm or}
\ \ \ \ \ \ \ \ \ \ \ \ \ \ \ \ \ \ \ \ \
\mathcal{C}^\infty,
\ \ \ \ \ \ \ \ \ \ \ \ \ \ \ \ \ \ \ \ \
\text{\rm or}
\ \ \ \ \ \ \ \ \ \ \ \ \ \ \ \ \ \ \ \ \
\mathcal{C}^\omega.
\]

Pick a point:
\[
p\in M,
\]
and take a (small) open neighborhood:
\[
p
\in
{\sf U}_p
\subset
\C^3.
\]

By definition (\cite{ Merker-Pocchiola-Sabzevari-5-CR-II,
Merker-5-CR-III}):
\[
\Big(
M^4
\,\subset\,
\C^3
\Big)
\,\,\in\,\,
\text{\sf General Class $\text{\sf II}$},
\]
if:
\[
\C\otimes_\R TM
=
T^{1,0}M
+
T^{0,1}M
+
\big[
T^{1,0}M,\,T^{0,1}M\big]
+
\big[T^{1,0}M,\,
\big[
T^{1,0}M,\,T^{0,1}M\big]\big].
\]

This means that for any local vector field generator:
\[
\mathcal{L}
\,=\,
\text{\rm section of}\,\,
T^{1,0}\big(M\cap{\sf U}_p\big),
\]
one has at every point $q \in M\cap {\sf U}_p$:
\[
{\bf 4}
=
\rank_\C
\Big(
\mathcal{L}\big\vert_q,\,\,
\overline{\mathcal{L}}\big\vert_q,\,\,
\big[\mathcal{L},\overline{\mathcal{L}}\big]
\Big\vert_q,\,\,
\big[\mathcal{L},\,\big[\mathcal{L},\overline{\mathcal{L}}\big]\big]
\Big\vert_q
\Big).
\]

\medskip

Next, take any (local) biholomorphism:
\[
h\colon\ \ \
{\sf U}_p
\overset{\sim}{\,\longrightarrow\,}
{\sf U}_{p'}'
=
h({\sf U}_p)
\ \ \ \ \ \ \ \ \ \ \ \ \
{\scriptstyle{(p'\,=\,\,h(p))}},
\]
which, when ${\sf U}_p$ is small enough, certainly 
transfers $M \cap {\sf U}_p$ to a certain CR-generic
submanifold:
\[
{M'}^4
\,:=\,
h\big(M\cap{\sf U}_p\big)
\,\subset\,
{\C'}^3,
\]
having of course the same smoothness:
\[
\mathcal{C}^\kappa\ \
(\kappa\geqslant 3),
\ \ \ \ \ \ \ \ \ \ \ \ \ \ \ \ \ \ \ \ \
\text{\rm or}
\ \ \ \ \ \ \ \ \ \ \ \ \ \ \ \ \ \ \ \ \
\mathcal{C}^\infty,
\ \ \ \ \ \ \ \ \ \ \ \ \ \ \ \ \ \ \ \ \
\text{\rm or}
\ \ \ \ \ \ \ \ \ \ \ \ \ \ \ \ \ \ \ \ \
\mathcal{C}^\omega.
\]

\medskip

Take also any local vector field generator:
\[
\mathcal{L}'
\,=\,
\text{\rm section of}\,\,
T^{1,0}M'.
\]
Then necessarily at every point $q' \in M'$, one 
also has (exercise, or {\em see} below):
\[
{\bf 4}
=
\rank_\C
\Big(
\mathcal{L}'\big\vert_{q'},\,\,
\overline{\mathcal{L}}'\big\vert_{q'},\,\,
\big[\mathcal{L}',\overline{\mathcal{L}}'\big]
\Big\vert_{q'},\,\,
\big[\mathcal{L}',\,\big[\mathcal{L}',\overline{\mathcal{L}}'\big]\big]
\Big\vert_{q'}
\Big).
\]

Write the components of $h$ and the target coordinates as:
\[
\aligned
(z,w_1,w_2)
&
\,\longmapsto\,
\big(z'(z,w_1,w_2),\,w_1'(z,w_1,w_2),\,w_2'(z,w_1,w_2)\big)
\\
&\ \ \ \
=:
(z',w_1',w_2').
\endaligned
\]

Consider the {\em inverse} of $h$:
\[
h'\colon \ \ \ \ \ \ \ \ \ \ \
\aligned
{\sf U}_{p'}'
&
\overset{\sim}{\,\longrightarrow\,}
{\sf U}_p
\ \ \ \ \ \ \ \ \ \ \ \ \ \ \ \ \ \ \ \ \ \ \ \ \ \ \ \ \ \ \ \ \ \ \
\ \ \ \ \ \ \ \ \ \ \ \ \ \ \ \ \ \ \ \
{\scriptstyle{(p\,=\,\,h'(p'))}},
\\
(z',w_1',w_2')
&
\,\longmapsto\,
\big(z(z',w_1',w_2'),\,w_1(z',w_1',w_2'),\,w_2(z',w_1',w_2')\big)
\\
&
\ \ \ \ \,
=
(z,w_1,w_2).
\endaligned
\]

Then there
exists a nowhere vanishing function:
\[
a\colon\ \ \
M\cap{\sf U}_p
\,\longrightarrow\,
\C\backslash\{0\},
\]
such that:
\[
h_*'
\big(\mathcal{L}'\big)
=
a\,\mathcal{L}.
\]

Simultaneously:
\[
h_*'
\big(\overline{\mathcal{L}}'\big)
=
\overline{a}\,\overline{\mathcal{L}}.
\]

Now setting:
\[
\aligned
\mathcal{T}
&
:=
\isqrt\,
\big[\mathcal{L},\,\overline{\mathcal{L}}\big],
\\
\mathcal{T}'
&
:=
\isqrt\,
\big[\mathcal{L}',\,\overline{\mathcal{L}}'\big],
\endaligned
\]
and setting:
\[
\aligned
\mathcal{S}
&
:=
\big[\mathcal{L},\mathcal{T}\big],
\\
\mathcal{S}'
&
:=
\big[\mathcal{L}',\mathcal{T}'\big],
\endaligned
\]
one has two (local) frames:
\[
\aligned
&
\big\{\mathcal{L},\,\overline{\mathcal{L}},\,
\mathcal{T},\,\mathcal{S}\big\}
\ \ \ \ \ \ \ \,
\text{\rm for}\ \
\C\otimes_\R T(M\cap{\sf U}_p),
\\
&
\big\{\mathcal{L}',\,\overline{\mathcal{L}}',\,
\mathcal{T}',\mathcal{S}'\big\}
\ \ \ \ \
\text{\rm for}\ \
\C\otimes_\R TM'.
\endaligned
\]

The Lie bracket $\mathcal{ T}'$ transfers
back to $M \cap {\sf U}_p$ through $h_*'$ as:
\[
\aligned
h_*'\big(\mathcal{T}'\big)
&
=
h_*'\big(\isqrt\,\big[\mathcal{L}',\,\overline{\mathcal{L}}'\big]\big)
\\
&
=
\isqrt\,
\big[
h_*'\big(\mathcal{L}'\big),\,
h_*'\big(\overline{\mathcal{L}}'\big)
\big]
\\
&
=
\isqrt\,
\big[a\,\mathcal{L},\,\,
\overline{a}\,\overline{\mathcal{L}}\big]
\\
&
=
a\overline{a}\,
\underbrace{\isqrt\,\,\big[\mathcal{L},\,\overline{\mathcal{L}}\big]}_{
=\,\mathcal{T}}
+
\underbrace{
\isqrt\,a\,\mathcal{L}\big(\overline{a}\big)}_{
=:\,\overline{b}}
\cdot
\overline{\mathcal{L}}
\underbrace{-
\isqrt\,\overline{a}\,\overline{\mathcal{L}}(a)}_{
=:\,b}
\cdot
\mathcal{L},
\endaligned
\]
so that, if one decides to set:
\[
b
:=
-\,\isqrt\,\overline{a}\,\overline{\mathcal{L}}(a),
\]
forgetting how this coefficient-function
is related to $a$, one obtains, exactly as for 
the General Class $\text{\sf I}$:
\[
h_*'\big(\mathcal{T}'\big)
=
a\overline{a}\,\mathcal{T}
+
\overline{b}\,\overline{\mathcal{L}}
+
b\,\mathcal{L}.
\]

\medskip

Next:
\[
\aligned
h_*'\big(\mathcal{S}'\big)
&
=
h_*'\big(\big[\mathcal{L}',\mathcal{T}'\big]\big)
\\
&
=
\big[
h_*'\big(\mathcal{L}'\big),\,
h_*'\big(\mathcal{T}'\big)
\big]
\\
&
=
\big[a\,\mathcal{L},\,\,
a\overline{a}\,\mathcal{T}
+
\overline{b}\,\overline{\mathcal{L}}
+
b\,\mathcal{L}
\big]
\\
&
=
aa\overline{a}\,
\underbrace{\big[\mathcal{L},\mathcal{T}\big]}_{
=\,\mathcal{S}}
+
a\overline{b}\,
\underbrace{\big[\mathcal{L},\overline{\mathcal{L}}\big]}_{
=\,-\,\isqrt\,\mathcal{T}}
+
ab\,\zero{\big[\mathcal{L},\mathcal{L}\big]}
+
\\
&
\ \ \ \ \
+
a\,\mathcal{L}\big(a\overline{a}\big)
\cdot
\mathcal{T}
+
a\,\mathcal{L}\big(\overline{b}\big)
\cdot
\overline{\mathcal{L}}
+
a\,\mathcal{L}(b)\cdot\mathcal{L}
-
\\
&
\ \ \ \ \
-\,
a\overline{a}\,\mathcal{T}(a)
\cdot
\mathcal{L}
-
\overline{b}\,\overline{\mathcal{L}}(a)
\cdot
\mathcal{L}
-
b\,\mathcal{L}(a)
\cdot
\mathcal{L},
\endaligned
\]
which is:
\[
\aligned
h_*'\big(\mathcal{S}'\big)
&
=
aa\overline{a}\cdot\mathcal{S}
+
\big(
-\,\isqrt\,a\overline{b}
+
a\,\mathcal{L}\big(a\overline{a}\big)
\big)
\cdot
\mathcal{T}
+
\big(a\,\mathcal{L}(\overline{b})\big)
\cdot
\overline{\mathcal{L}}
+
\\
&
\ \ \ \ \
+
\big(
a\,\mathcal{L}(b)
-
a\overline{a}\,\mathcal{T}(a)
-
\overline{b}\,\overline{\mathcal{L}}(a)
-
b\,\mathcal{L}(a)
\big)
\cdot
\mathcal{L},
\endaligned
\]
so that if one decides to set:
\[
\aligned
c
&
:=
-\,\isqrt\,a\overline{b}
+
a\,\mathcal{L}\big(a\overline{a}\big),
\\
d
&
:=
a\,\mathcal{L}\big(\overline{b}\big),
\\
e
&
:=
a\,\mathcal{L}(b)
-
a\overline{a}\,\mathcal{T}(a)
-
\overline{b}\,\overline{\mathcal{L}}(a)
-
b\,\mathcal{L}(a),
\endaligned
\]
one has in:

\medskip\noindent{\bf Summary.}
{\em Through any local biholomorphic equivalences
between CR-generic submanifolds of $\C^3$ belonging
to the General Class} $\text{\sf II}$:
\[
M^4
\overset{\sim}{\,\longrightarrow\,}
{M'}^4,
\]
{\em for any two choices of local vector
field generators:}
\[
\aligned
&
\mathcal{L}
\ \ \ \ \
\text{\rm for}\ \
T^{1,0}M,
\\
&
\mathcal{L}'
\ \ \ \ \
\text{\rm for}\ \
T^{1,0}M',
\endaligned
\]
{\em the transfer of frame obeys the rule:}
\[
\left(\!\!
\begin{array}{c}
\mathcal{L}'
\\
\overline{\mathcal{L}}'
\\
\mathcal{T}'
\\
\mathcal{S}'
\end{array}
\!\!\right)
=
\left(\!
\begin{array}{cccc}
a & 0 & 0 & 0
\\
0 & \overline{a} & 0 & 0 
\\
b & \overline{b} & a\overline{a} & 0 
\\
e & d & c & aa\overline{a}
\end{array}
\!\right)
\left(\!\!
\begin{array}{c}
\mathcal{L}
\\
\overline{\mathcal{L}}
\\
\mathcal{T}
\\
\mathcal{S}
\end{array}
\!\!\right),
\]
{\em for some five local functions:}
\[
\aligned
a\colon\ \ \
M
&
\,\longrightarrow\,
\C\backslash\{0\},
\\
b,c,d,e\colon\ \ \
M
&
\,\longrightarrow\,
\C.
\endaligned
\]

\medskip

This means that the {\em ambiguity} in the choice of a local frame:
\[
\big\{
\mathcal{L},\,
\overline{\mathcal{L}},\,
\mathcal{T},\,
\mathcal{S}
\big\}
\ \ \ \ \
\text{\rm for}\ \
\C\otimes_\R TM
\]
which comes\,\,---\,\,{\em naturally from the
point of view of CR geometry}\,\,---\,\,from a choice of a local frame:
\[
\mathcal{L}
\ \ \ \ \
\text{\rm for}\ \
T^{1,0}M
\]
is represented by general changes of frames whose
matrices are of the form: 
\[
\left(\!
\begin{array}{cccc}
a & 0 & 0 & 0
\\
0 & \overline{a} & 0 & 0
\\
b & \overline{b} & a\overline{a} & 0
\\
e & d & c & aa\overline{a}
\end{array}
\!\right),
\]
the entries being coefficient-functions depending on existing
equivalences $M \overset{ \sim}{ \longrightarrow} M'$,
and\big/or on the choice of local coordinates.

\medskip

Thus, the {\em initial
$G$-structure} for the biholomorphic equivalence problem between
CR-generic submanifolds $M^4 \subset \C^3$ belonging to the General Class
$\text{\sf II}$ is a reduction of the full linear group 
${\sf GL}_4 ( \C)$ to
the mentioned subgroup in which the\,\,---\,\,possibly
unknown\,\,---\,\,functions $a$, $b$,
$c$, $d$, $e$ are {\em replaced} by independent
complex variables.

\medskip\noindent{\bf Lemma.}
{\em The set of matrices:}
\[
{\sf G}_{\text{\sf II}}^{\sf initial}
\,:=\,
\left\{
\left(\!\!
\begin{array}{cccc}
{\sf a} & 0 & 0 & 0
\\
0 & \overline{\sf a} & 0 & 0
\\
{\sf b} & \overline{\sf b} & {\sf a}\overline{\sf a} & 0
\\
{\sf e} & {\sf d} & {\sf c} & {\sf a}{\sf a}\overline{\sf a}
\end{array}
\!\!\right)
\,\in\,
\mathcal{M}_{4\times 4}(\C)\,\colon\,\,
{\sf a}\,\in\,\C\backslash\{0\},\,\,
{\sf b},\,{\sf c},\,{\sf d},\,{\sf e}\,\in\,\C
\right\}
\]
{\em is a (closed) $10$-dimensional
real matrix subgroup of the full:}
\[
{\sf GL}_4(\C)
\,:=\,
\left\{
\pi
=
\left(\!\!
\begin{array}{cccc}
\pi_{1,1} & \pi_{1,2} & \pi_{1,3} & \pi_{1,4}
\\
\pi_{2,1} & \pi_{2,2} & \pi_{3,3} & \pi_{2,4}
\\
\pi_{3,1} & \pi_{3,2} & \pi_{3,3} & \pi_{3,4}
\\
\pi_{4,1} & \pi_{4,2} & \pi_{4,3} & \pi_{4,4}
\end{array}
\!\!\right)
\,\in\,
\mathcal{M}_{4\times 4}(\C)\,\colon\,\,
0
\neq
\det\,\pi
\right\}.
\]

\proof
Closedness under multiplication (composition):
\[
\aligned
&
\left(\!\!
\begin{array}{cccc}
{\sf a}_1 & 0 & 0 & 0
\\
0 & \overline{\sf a}_1 & 0 & 0
\\
{\sf b}_1 & \overline{\sf b}_1 &
{\sf a}_1\overline{\sf a}_1 & 0
\\
{\sf e}_1 & {\sf d}_1 & {\sf c}_1 & {\sf a}_1{\sf a}_1
\overline{\sf a}_1
\end{array}
\!\!\right)
\cdot
\left(\!\!
\begin{array}{cccc}
{\sf a}_2 & 0 & 0 & 0
\\
0 & \overline{\sf a}_2 & 0 & 0
\\
{\sf b}_2 & \overline{\sf b}_2 &
{\sf a}_2\overline{\sf a}_2 & 0
\\
{\sf e}_2 & {\sf d}_2 & {\sf c}_2 & {\sf a}_2{\sf a}_2
\overline{\sf a}_2
\end{array}
\!\!\right)
\\
&
\ \ \ \ \ \ \ \ \ \
=
\left(\!\!
\begin{array}{cccc}
{\sf a}_1{\sf a}_2 & 0 & 0 & 0
\\
0 & \overline{\sf a}_1\overline{\sf a}_2 & 0 & 0
\\
{\sf b}_1{\sf a}_2
\!+\!{\sf a}_1\overline{\sf a}_1{\sf b}_2 & 
\overline{\sf b}_1\overline{\sf a}_2
\!+\!{\sf a}_1\overline{\sf a}_1\overline{\sf b}_2 &
{\sf a}_1\overline{\sf a}_1{\sf a}_2\overline{\sf a}_2 & 0\medskip
\\
\substack{{\sf e}_1{\sf a}_2+{\sf c}_1{\sf b}_2+\\
+{\sf a}_1{\sf a}_1\overline{\sf a}_1{\sf e}_2}
&
\substack{{\sf d}_1\overline{\sf a}_2+{\sf c}_1\overline{\sf b}_2+\\
+{\sf a}_1{\sf a}_1\overline{\sf a}_1{\sf d}_2}
&
\substack{{\sf c}_1{\sf a}_2\overline{\sf a}_2+\\
{\sf a}_1{\sf a}_1\overline{\sf a}_1{\sf c}_2}
&
\substack{
{\sf a}_1{\sf a}_1\overline{\sf a}_1
{\sf a}_2{\sf a}_2\overline{\sf a}_2}
\end{array}
\!\!\right)
\\
&
\ \ \ \ \ \ \ \ \ \
=:
\left(\!\!
\begin{array}{cccc}
{\sf a}_3 & 0 & 0 & 0
\\
0 & \overline{\sf a}_3 & 0 & 0
\\
{\sf b}_3 & \overline{\sf b}_3 &
{\sf a}_3\overline{\sf a}_3 & 0
\\
{\sf e}_3 & {\sf d}_3 & {\sf c}_3 & {\sf a}_3{\sf a}_3\overline{\sf a}_3
\end{array}
\!\!\right)
\endaligned
\]
is visibly clear, after setting:
\[
\aligned
{\sf a}_3
&
:=
{\sf a}_1{\sf a}_2
\ \ \ \ \ \ \ \ \ \ \ \ \ \ \ \ \ \ \ \ \ \ \ \ \
{\scriptstyle{(\neq\,0,\,\,
\text{\rm again})}},
\\
{\sf b}_3
&
:=
{\sf b}_1{\sf a}_2
+
{\sf a}_1\overline{\sf a}_1{\sf b}_2,
\\
{\sf c}_3
&
:=
{\sf c}_1{\sf a}_2\overline{\sf a}_2
+
{\sf a}_1{\sf a}_1\overline{\sf a}_1{\sf c}_2,
\\
{\sf d}_3
&
:=
{\sf d}_1\overline{\sf a}_2
+
{\sf c}_1\overline{\sf b}_2
+
{\sf a}_1{\sf a}_1\overline{\sf a}_1{\sf d}_2,
\\
{\sf e}_3
&
:=
{\sf e}_1{\sf a}_2
+
{\sf c}_1{\sf b}_2
+
{\sf a}_1{\sf a}_1\overline{\sf a}_1{\sf e}_2.
\endaligned
\]

Quite similarly, the inverse:
\[
\aligned
\left(\!\!
\begin{array}{cccc}
{\sf a} & 0 & 0 & 0
\\
0 & \overline{\sf a} & 0 & 0
\\
{\sf b} & \overline{\sf b} &
{\sf a}\overline{\sf a} & 0
\\
{\sf e} & {\sf d} & {\sf c} & 
{\sf a}{\sf a}\overline{\sf a}
\end{array}
\!\!\right)^{-1}
&
=
\left(\!\!
\begin{array}{cccc}
\frac{1}{\sf a} & 0 & 0 & 0
\\
0 & \frac{1}{\overline{\sf a}} & 0 & 0\smallskip
\\
\frac{-{\sf b}}{{\sf a}{\sf a}\overline{\sf a}} & 
\frac{-\overline{\sf b}}{{\sf a}\overline{\sf a}\overline{\sf a}} &
\frac{1}{{\sf a}\overline{\sf a}} & 0\smallskip
\\
\frac{{\sf b}{\sf c}}{{\sf a}^4\overline{\sf a}^2}
-
\frac{{\sf e}}{{\sf a}^3\overline{\sf a}}
&
\frac{{\sf c}\overline{\sf b}}{{\sf a}^3
\overline{\sf a}^3}
-
\frac{{\sf d}}{{\sf a}^2\overline{\sf a}^2}
&
\frac{-{\sf c}}{{\sf a}^3\overline{\sf a}^2}
&
\frac{1}{{\sf a}^2\overline{\sf a}}
\end{array}
\!\!\right)
\\
&
=:
\left(\!\!
\begin{array}{cccc}
{\sf a}^\sim & 0 & 0 & 0
\\
0 & \overline{\sf a}^\sim & 0 & 0
\\
{\sf b}^\sim & \overline{\sf b}^\sim &
{\sf a}^\sim\overline{\sf a}^\sim & 0
\\
{\sf e}^\sim & {\sf d}^\sim & {\sf c}^\sim &
{\sf a}^\sim{\sf a}^\sim\overline{\sf a}^\sim
\end{array}
\!\!\right),
\endaligned
\]
also belongs to the subgroup, with:
\[
\aligned
{\sf a}^\sim
&
:=
\frac{1}{\sf a}
\ \ \ \ \ \ \ \ \ \ \ \ \ \ \ \ \ \ \ \ \ \ \ \ \
{\scriptstyle{(\neq\,0)}},
\\
{\sf b}^\sim
&
:=
-\,\frac{\sf b}{{\sf a}{\sf a}\overline{\sf a}},
\\
{\sf c}^\sim
&
:=
-\,\frac{{\sf c}}{{\sf a}{\sf a}{\sf a}\overline{\sf a}\overline{\sf a}},
\\
{\sf d}^\sim
&
:=
\frac{{\sf c}\overline{\sf b}}{
{\sf a}{\sf a}{\sf a}\overline{\sf a}\overline{\sf a}\overline{\sf a}}
-
\frac{{\sf d}}{{\sf a}{\sf a}\overline{\sf a}\overline{\sf a}},
\\
{\sf e}^\sim
&
:=
\frac{{\sf b}{\sf c}}{{\sf a}{\sf a}{\sf a}{\sf a}
\overline{\sf a}\overline{\sf a}}
-
\frac{{\sf e}}{{\sf a}{\sf a}{\sf a}\overline{\sf a}},
\endaligned
\]
which concludes.
\endproof

\noindent{\bf Proposition.}
{\em On a $4$-dimensional CR-generic submanifold:}
\[
\Big(
M^4
\,\subset\,
\C^3
\Big)
\,\,\in\,\,
\text{\sf General Class $\text{\sf II}$},
\]
{\em having biholomorphically invariant $(1, 0)$ CR bundle:}
\[
T^{1,0}M
\,\subset\,
\C\otimes_\R TM,
\]
{\em for any choice of local vector field generator:}
\[
\mathcal{L}
\ \ \ \ \
\text{\rm for}\ \
T^{1,0}M,
\]
{\em the associated frame:}
\[
\big\{
\mathcal{L},\,\overline{\mathcal{L}},\,
\isqrt\,\big[\mathcal{L},\overline{\mathcal{L}}\big],\,
\big[\mathcal{L},\,\isqrt\big[\mathcal{L},\overline{\mathcal{L}}\big]\big]
\big\}
\,=:\,
\big\{
\mathcal{L},\,\overline{\mathcal{L}},\,\mathcal{T},\,\mathcal{S}
\big\}
\]
{\em for the full complexified tangent bundle:}
\[
\C\otimes_\R TM
\]
{\em performs a reduction of the full ${\sf GL}_4 ( \C)$-structure
of $\C \otimes_\R TM$ to the $10$-dimensional subgroup:}
\[
{\sf G}_{\text{\sf II}}^{\sf initial}
\,:=\,
\left\{
\left(\!\!
\begin{array}{cccc}
{\sf a} & 0 & 0 & 0
\\
0 & \overline{\sf a} & 0 & 0
\\
{\sf b} & \overline{\sf b} & {\sf a}\overline{\sf a} & 0
\\
{\sf e} & {\sf d} & {\sf c} & {\sf a}{\sf a}\overline{\sf a}
\end{array}
\!\!\right)
\,\in\,
\mathcal{M}_{4\times 4}(\C)\,\colon\,\,
{\sf a}\,\in\,\C\backslash\{0\},\,\,
{\sf b},\,{\sf c},\,{\sf d},\,{\sf e}\,\in\,\C
\right\}.
\qed
\]


\bigskip

\section{\sf General class $\text{\sf III}_{\text{\sf 1}}$}
\label{general-class-III-1}
\HEAD{\ref{general-class-III-1}.~General class 
$\text{\sf III}_{\text{\sf 1}}$}{
Jo\"el {\sc Merker}, D\'epartement de Math\'ematiques d'Orsay}

\medskip

Equip $\C^4$ with coordinates:
\[
(z,w_1,w_2,w_3)
\,\in\,
\C^4.
\]

Let a connected CR-generic submanifold:
\[
M^5
\subset
\C^4
\]
be of smoothness:
\[
\mathcal{C}^\kappa\ \
(\kappa\geqslant 3),
\ \ \ \ \ \ \ \ \ \ \ \ \ \ \ \ \ \ \ \ \
\text{\rm or}
\ \ \ \ \ \ \ \ \ \ \ \ \ \ \ \ \ \ \ \ \
\mathcal{C}^\infty,
\ \ \ \ \ \ \ \ \ \ \ \ \ \ \ \ \ \ \ \ \
\text{\rm or}
\ \ \ \ \ \ \ \ \ \ \ \ \ \ \ \ \ \ \ \ \
\mathcal{C}^\omega.
\]

Pick a point:
\[
p\in M,
\]
and take a (small) open neighborhood:
\[
p
\in
{\sf U}_p
\subset
\C^4.
\]

By definition (\cite{ Merker-Pocchiola-Sabzevari-5-CR-II,
Merker-5-CR-III}):
\[
\Big(
M^5
\,\subset\,
\C^4
\Big)
\,\,\in\,\,
\text{\sf General Class $\text{\sf III}_{\text{\sf 1}}$},
\]
if:
\[
\aligned
\C\otimes_\R TM
=
T^{1,0}M
+
T^{0,1}M
+
\big[
T^{1,0}M,\,T^{0,1}M\big]
&
+
\big[T^{1,0}M,\,
\big[
T^{1,0}M,\,T^{0,1}M\big]\big]
\\
&
+
\big[T^{0,1}M,\,
\big[
T^{1,0}M,\,T^{0,1}M\big]\big].
\endaligned
\]

This means that for any local vector field generator:
\[
\mathcal{L}
\,=\,
\text{\rm section of}\,\,
T^{1,0}\big(M\cap{\sf U}_p\big),
\]
one has at every point $q \in M\cap {\sf U}_p$:
\[
{\bf 5}
=
\rank_\C
\Big(
\mathcal{L}\big\vert_q,\,\,
\overline{\mathcal{L}}\big\vert_q,\,\,
\big[\mathcal{L},\overline{\mathcal{L}}\big]
\Big\vert_q,\,\,
\big[\mathcal{L},\,\big[\mathcal{L},\overline{\mathcal{L}}\big]\big]
\Big\vert_q,\,\,
\big[\overline{\mathcal{L}},\,\big[\mathcal{L},
\overline{\mathcal{L}}\big]\big]
\Big\vert_q
\Big).
\]

\medskip

Next, take any (local) biholomorphism:
\[
h\colon\ \ \
{\sf U}_p
\overset{\sim}{\,\longrightarrow\,}
{\sf U}_{p'}'
=
h({\sf U}_p)
\ \ \ \ \ \ \ \ \ \ \ \ \
{\scriptstyle{(p'\,=\,\,h(p))}},
\]
which, when ${\sf U}_p$ is small enough, certainly 
transfers $M \cap {\sf U}_p$ to a certain CR-generic
submanifold:
\[
{M'}^5
\,:=\,
h\big(M\cap{\sf U}_p\big)
\,\subset\,
{\C'}^4,
\]
having of course the same smoothness:
\[
\mathcal{C}^\kappa\ \
(\kappa\geqslant 3),
\ \ \ \ \ \ \ \ \ \ \ \ \ \ \ \ \ \ \ \ \
\text{\rm or}
\ \ \ \ \ \ \ \ \ \ \ \ \ \ \ \ \ \ \ \ \
\mathcal{C}^\infty,
\ \ \ \ \ \ \ \ \ \ \ \ \ \ \ \ \ \ \ \ \
\text{\rm or}
\ \ \ \ \ \ \ \ \ \ \ \ \ \ \ \ \ \ \ \ \
\mathcal{C}^\omega.
\]

\medskip

Take also any local vector field generator:
\[
\mathcal{L}'
\,=\,
\text{\rm section of}\,\,
T^{1,0}M'.
\]
Then necessarily at every point $q' \in M'$, one 
also has (exercise, or {\em see} below):
\[
{\bf 5}
=
\rank_\C
\Big(
\mathcal{L}'\big\vert_{q'},\,\,
\overline{\mathcal{L}}'\big\vert_{q'},\,\,
\big[\mathcal{L}',\overline{\mathcal{L}}'\big]
\Big\vert_{q'},\,\,
\big[\mathcal{L}',\,\big[\mathcal{L}',\overline{\mathcal{L}}'\big]\big]
\Big\vert_{q'},\,\,
\big[\overline{\mathcal{L}}',\,
\big[\mathcal{L}',\overline{\mathcal{L}}'\big]\big]
\Big\vert_{q'}
\Big).
\]

Write the components of $h$ and the target coordinates as:
\[
\!\!\!\!\!\!\!\!\!\!\!\!\!\!\!\!\!\!\!\!
\aligned
(z,w_1,w_2,w_3)
&
\,\longmapsto\,
\big(z'(z,w_1,w_2,w_3),\,w_1'(z,w_1,w_2,w_3),\,w_2'(z,w_1,w_2,w_3),\,
w_3'(z,w_1,w_2,w_3)\big)
\\
&\ \ \ \
=:
(z',w_1',w_2',w_3').
\endaligned
\]

Consider the {\em inverse} $h'$
of $h$:
\[
\!\!\!\!\!\!\!\!\!\!\!\!\!\!\!\!\!\!\!\!
\aligned
{\sf U}_{p'}'
&
\overset{\sim}{\,\longrightarrow\,}
{\sf U}_p
\ \ \ \ \ \ \ \ \ \ \ \ \ \ \ \ \ \ \ \ \ \ \ \ \ \ \ \ \ \ \ \ \ \ \
\ \ \ \ \ \ \ \ \ \ \ \ \ \ \ \ \ \ \ \ \ \ \ \ \ \ \ \ \ \ \ \ \ \ \
\ \ \ \ \ \ \ \ \ \ \ \ \ \ \ \ \ \ \ \ \ \ \ \ \ \ \ \ \ \ \ \ \ 
{\scriptstyle{(p\,=\,\,h'(p'))}},
\\
(z',w_1',w_2',w_3')
&
\,\longmapsto\,
\big(z(z',w_1',w_2',w_3'),\,w_1(z',w_1',w_2',w_3'),\,
w_2(z',w_1',w_2',w_3'),\,w_3(z',w_1',w_2',w_3')\big)
\\
&\ \ \ \ \,
=
(z,w_1,w_2,w_3).
\endaligned
\]

Then there
exists a nowhere vanishing function:
\[
a\colon\ \ \
M\cap{\sf U}_p
\,\longrightarrow\,
\C\backslash\{0\},
\]
such that:
\[
h_*'
\big(\mathcal{L}'\big)
=
a\,\mathcal{L}.
\]

Simultaneously:
\[
h_*'
\big(\overline{\mathcal{L}}'\big)
=
\overline{a}\,\overline{\mathcal{L}}.
\]

Now setting:
\[
\aligned
\mathcal{T}
&
:=
\isqrt\,
\big[\mathcal{L},\,\overline{\mathcal{L}}\big],
\\
\mathcal{T}'
&
:=
\isqrt\,
\big[\mathcal{L}',\,\overline{\mathcal{L}}'\big],
\endaligned
\]
setting:
\[
\aligned
\mathcal{S}
&
:=
\big[\mathcal{L},\mathcal{T}\big],
\\
\mathcal{S}'
&
:=
\big[\mathcal{L}',\mathcal{T}'\big],
\endaligned
\]
whence:
\[
\aligned
\overline{\mathcal{S}}
&
=
\big[\overline{\mathcal{L}},\,\mathcal{T}\big],
\\
\overline{\mathcal{S}}'
&
=
\big[\overline{\mathcal{L}}',\,\mathcal{T}'\big],
\endaligned
\]
one has two (local) frames:
\[
\aligned
&
\big\{\mathcal{L},\,\overline{\mathcal{L}},\,
\mathcal{T},\,\mathcal{S},\,\overline{\mathcal{S}}\big\}
\ \ \ \ \ \ \ \ \ \,
\text{\rm for}\ \
\C\otimes_\R T(M\cap{\sf U}_p),
\\
&
\big\{\mathcal{L}',\,\overline{\mathcal{L}}',\,
\mathcal{T}',\,\mathcal{S}',\,\overline{\mathcal{S}}'\big\}
\ \ \ \ \
\text{\rm for}\ \
\C\otimes_\R TM'.
\endaligned
\]

As for the General Class $\text{\sf II}$:
\[
\aligned
h_*'\big(\mathcal{T}'\big)
&
=
a\overline{a}\,\mathcal{T}
+
\overline{b}\,\overline{\mathcal{L}}
+
b\,\mathcal{L},
\\
h_*'\big(\mathcal{S}'\big)
&
=
aa\overline{a}\,\mathcal{S}
+
c\,\mathcal{T}
+
d\,\overline{\mathcal{L}}
+
e\,\mathcal{L},
\endaligned
\]
whence by plain conjugation:
\[
h_*'\big(\overline{\mathcal{S}}'\big)
=
a\overline{a}\overline{a}\,\overline{\mathcal{S}}
+
\overline{c}\,\mathcal{T}
+
\overline{e}\,\overline{\mathcal{L}}
+
\overline{d}\,\mathcal{L}.
\]

\medskip\noindent{\bf Summary.}
{\em Through any local biholomorphic equivalences
between CR-generic submanifolds of $\C^4$ belonging
to the General Class} $\text{\sf III}_{\text{\sf 1}}$:
\[
M^5
\overset{\sim}{\,\longrightarrow\,}
{M'}^5,
\]
{\em for any two choices of local vector
field generators:}
\[
\aligned
&
\mathcal{L}
\ \ \ \ \
\text{\rm for}\ \
T^{1,0}M,
\\
&
\mathcal{L}'
\ \ \ \ \
\text{\rm for}\ \
T^{1,0}M',
\endaligned
\]
{\em the transfer of frame obeys the rule:}
\[
\left(\!\!
\begin{array}{c}
\mathcal{L}'
\\
\overline{\mathcal{L}}'
\\
\mathcal{T}'
\\
\mathcal{S}'
\\
\overline{\mathcal{S}}'
\end{array}
\!\!\right)
=
\left(\!
\begin{array}{ccccc}
a & 0 & 0 & 0 & 0
\\
0 & \overline{a} & 0 & 0 & 0
\\
b & \overline{b} & a\overline{a} & 0 & 0
\\
e & d & c & aa\overline{a} & 0
\\
\overline{d} & \overline{e} & \overline{c} & 0 &
a\overline{a}\overline{a}
\end{array}
\!\right)
\left(\!\!
\begin{array}{c}
\mathcal{L}
\\
\overline{\mathcal{L}}
\\
\mathcal{T}
\\
\mathcal{S}
\\
\overline{\mathcal{S}}
\end{array}
\!\!\right),
\]
{\em for some five local functions:}
\[
\aligned
a\colon\ \ \
M
&
\,\longrightarrow\,
\C\backslash\{0\},
\\
b,c,d,e\colon\ \ \
M
&
\,\longrightarrow\,
\C.
\endaligned
\]

\medskip

This means that the {\em ambiguity} in the choice of a local frame:
\[
\big\{
\mathcal{L},\,
\overline{\mathcal{L}},\,
\mathcal{T},\,
\mathcal{S},\,
\overline{\mathcal{S}}
\big\}
\ \ \ \ \
\text{\rm for}\ \
\C\otimes_\R TM
\]
which comes\,\,---\,\,{\em naturally from the
point of view of CR geometry}\,\,---\,\,from a choice of a local frame:
\[
\mathcal{L}
\ \ \ \ \
\text{\rm for}\ \
T^{1,0}M
\]
is represented by general changes of frames whose
matrices are of the form: 
\[
\left(\!
\begin{array}{ccccc}
a & 0 & 0 & 0 & 0
\\
0 & \overline{a} & 0 & 0 & 0
\\
b & \overline{b} & a\overline{a} & 0 & 0
\\
e & d & c & aa\overline{a} & 0
\\
\overline{d} & \overline{e} & \overline{c} & 0 &
a\overline{a}\overline{a}
\end{array}
\!\right),
\]
the entries being coefficient-functions depending on existing
equivalences $M \overset{ \sim}{ \longrightarrow} M'$,
and\big/or on the choice of local coordinates.

\medskip

Thus, the {\em initial
$G$-structure} for the biholomorphic equivalence problem between
CR-generic submanifolds $M^5 \subset \C^4$ belonging to the General Class
$\text{\sf III}_{\text{\sf 1}}$ is a reduction of the full linear group 
${\sf GL}_5 ( \C)$ to
the mentioned subgroup in which the\,\,---\,\,possibly
unknown\,\,---\,\,functions $a$, $b$,
$c$, $d$, $e$ are {\em replaced} by independent
complex variables.

\medskip\noindent{\bf Lemma.}
{\em The set of matrices:}
\[
{\sf G}_{\text{\sf III}_{\text{\sf 1}}}^{\sf initial}
\,:=\,
\left\{
\left(\!\!
\begin{array}{ccccc}
{\sf a} & 0 & 0 & 0 & 0
\\
0 & \overline{\sf a} & 0 & 0 & 0
\\
{\sf b} & \overline{\sf b} & {\sf a}\overline{\sf a} & 0 & 0
\\
{\sf e} & {\sf d} & {\sf c} & {\sf a}{\sf a}\overline{\sf a} & 0
\\
\overline{\sf d} & \overline{\sf e} & \overline{\sf c} & 0 &
{\sf a}\overline{\sf a}\overline{\sf a}
\end{array}
\!\!\right)
\,\in\,
\mathcal{M}_{5\times 5}(\C)\,\colon\,\,
{\sf a}\,\in\,\C\backslash\{0\},\,\,
{\sf b},\,{\sf c},\,{\sf d},\,{\sf e}\,\in\,\C
\right\}
\]
{\em is a (closed) $10$-dimensional
real matrix subgroup of the full:}
\[
{\sf GL}_5(\C)
\,:=\,
\left\{
\pi
=
\left(\!\!
\begin{array}{ccccc}
\pi_{1,1} & \pi_{1,2} & \pi_{1,3} & \pi_{1,4} & \pi_{1,5}
\\
\pi_{2,1} & \pi_{2,2} & \pi_{3,3} & \pi_{2,4} & \pi_{2,5}
\\
\pi_{3,1} & \pi_{3,2} & \pi_{3,3} & \pi_{3,4} & \pi_{3,5}
\\
\pi_{4,1} & \pi_{4,2} & \pi_{4,3} & \pi_{4,4} & \pi_{4,5}
\\
\pi_{5,1} & \pi_{5,2} & \pi_{5,3} & \pi_{5,4} & \pi_{5,5}
\end{array}
\!\!\right)
\,\in\,
\mathcal{M}_{5\times 5}(\C)\,\colon\,\,
0
\neq
\det\,\pi
\right\}.
\]

\proof
Closedness under multiplication (composition):
\[
\aligned
&
\left(\!\!
\begin{array}{ccccc}
{\sf a}_1 & 0 & 0 & 0 & 0
\\
0 & \overline{\sf a}_1 & 0 & 0 & 0
\\
{\sf b}_1 & \overline{\sf b}_1 &
{\sf a}_1\overline{\sf a}_1 & 0 & 0
\\
{\sf e}_1 & {\sf d}_1 & {\sf c}_1 & {\sf a}_1{\sf a}_1
\overline{\sf a}_1 & 0
\\
\overline{\sf d}_1 & \overline{\sf e}_1 & \overline{\sf c}_1 & 0 &
{\sf a}_1{\sf a}_1\overline{\sf a}_1
\end{array}
\!\!\right)
\cdot
\left(\!\!
\begin{array}{ccccc}
{\sf a}_2 & 0 & 0 & 0 & 0
\\
0 & \overline{\sf a}_2 & 0 & 0 & 0
\\
{\sf b}_2 & \overline{\sf b}_2 &
{\sf a}_2\overline{\sf a}_2 & 0 & 0
\\
{\sf e}_2 & {\sf d}_2 & {\sf c}_2 & {\sf a}_2{\sf a}_2
\overline{\sf a}_2 & 0
\\
\overline{\sf d}_2 & \overline{\sf e}_2 & \overline{\sf c}_2 & 0 &
{\sf a}_2{\sf a}_2\overline{\sf a}_2
\end{array}
\!\!\right)
\\
&
\ \ \ \ \ \ \ \ \ \
=
\left(\!\!
\begin{array}{ccccc}
{\sf a}_1{\sf a}_2 & 0 & 0 & 0 & 0
\\
0 & \overline{\sf a}_1\overline{\sf a}_2 & 0 & 0 & 0
\\
{\sf b}_1{\sf a}_2
\!+\!{\sf a}_1\overline{\sf a}_1{\sf b}_2 & 
\overline{\sf b}_1\overline{\sf a}_2
\!+\!{\sf a}_1\overline{\sf a}_1\overline{\sf b}_2 &
{\sf a}_1\overline{\sf a}_1{\sf a}_2\overline{\sf a}_2 & 0 & 0\medskip
\\
\substack{{\sf e}_1{\sf a}_2+{\sf c}_1{\sf b}_2+\\
+{\sf a}_1{\sf a}_1\overline{\sf a}_1{\sf e}_2}
&
\substack{{\sf d}_1\overline{\sf a}_2+{\sf c}_1\overline{\sf b}_2+\\
+{\sf a}_1{\sf a}_1\overline{\sf a}_1{\sf d}_2}
&
\substack{{\sf c}_1{\sf a}_2\overline{\sf a}_2+\\
{\sf a}_1{\sf a}_1\overline{\sf a}_1{\sf c}_2}
&
\substack{
{\sf a}_1{\sf a}_1\overline{\sf a}_1
{\sf a}_2{\sf a}_2\overline{\sf a}_2} & 0\medskip
\\
\substack{
\overline{\sf d}_1{\sf a}_2+\overline{\sf c}_1{\sf b}_2+\\
+{\sf a}_1\overline{\sf a}_1\overline{\sf a}_1\overline{\sf d}_2}
&
\substack{
\overline{\sf e}_1\overline{\sf a}_2+\overline{\sf c}_1\overline{\sf b}_2+\\
+{\sf a}_1\overline{\sf a}_1\overline{\sf a}_1\overline{\sf e}_2}
&
\substack{
\overline{\sf c}_1{\sf a}_2\overline{\sf a}_2+\\
+{\sf a}_1\overline{\sf a}_1\overline{\sf a}_1\overline{\sf c}_2}
&
0
&
\substack{
{\sf a}_1\overline{\sf a}_1\overline{\sf a}_1
{\sf a}_2\overline{\sf a}_2\overline{\sf a}_2}
\end{array}
\!\!\right)
\\
&
\ \ \ \ \ \ \ \ \ \
=:
\left(\!\!
\begin{array}{ccccc}
{\sf a}_3 & 0 & 0 & 0 & 0
\\
0 & \overline{\sf a}_3 & 0 & 0 & 0
\\
{\sf b}_3 & \overline{\sf b}_3 &
{\sf a}_3\overline{\sf a}_3 & 0 & 0
\\
{\sf e}_3 & {\sf d}_3 & {\sf c}_3 & {\sf a}_3{\sf a}_3\overline{\sf a}_3 & 0
\\
\overline{\sf d}_3 & \overline{\sf e}_3 & \overline{\sf c}_3 & 0 &
{\sf a}_3\overline{\sf a}_3\overline{\sf a}_3
\end{array}
\!\!\right)
\endaligned
\]
is visibly clear, after setting:
\[
\aligned
{\sf a}_3
&
:=
{\sf a}_1{\sf a}_2
\ \ \ \ \ \ \ \ \ \ \ \ \ \ \ \ \ \ \ \ \ \ \ \ \
{\scriptstyle{(\neq\,0,\,\,
\text{\rm again})}},
\\
{\sf b}_3
&
:=
{\sf b}_1{\sf a}_2
+
{\sf a}_1\overline{\sf a}_1{\sf b}_2,
\\
{\sf c}_3
&
:=
{\sf c}_1{\sf a}_2\overline{\sf a}_2
+
{\sf a}_1{\sf a}_1\overline{\sf a}_1{\sf c}_2,
\\
{\sf d}_3
&
:=
{\sf d}_1\overline{\sf a}_2
+
{\sf c}_1\overline{\sf b}_2
+
{\sf a}_1{\sf a}_1\overline{\sf a}_1{\sf d}_2,
\\
{\sf e}_3
&
:=
{\sf e}_1{\sf a}_2
+
{\sf c}_1{\sf b}_2
+
{\sf a}_1{\sf a}_1\overline{\sf a}_1{\sf e}_2.
\endaligned
\]

Quite similarly, the inverse:
\[
\aligned
\left(\!\!
\begin{array}{ccccc}
{\sf a} & 0 & 0 & 0 & 0
\\
0 & \overline{\sf a} & 0 & 0 & 0
\\
{\sf b} & \overline{\sf b} &
{\sf a}\overline{\sf a} & 0 & 0
\\
{\sf e} & {\sf d} & {\sf c} & 
{\sf a}{\sf a}\overline{\sf a} & 0
\\
\overline{\sf d} & \overline{\sf e} & \overline{\sf c} & 0 &
{\sf a}\overline{\sf a}\overline{\sf a}
\end{array}
\!\!\right)^{-1}
&
=
\left(\!\!
\begin{array}{ccccc}
\frac{1}{\sf a} & 0 & 0 & 0 & 0
\\
0 & \frac{1}{\overline{\sf a}} & 0 & 0 & 0\smallskip
\\
\frac{-{\sf b}}{{\sf a}{\sf a}\overline{\sf a}} & 
\frac{-\overline{\sf b}}{{\sf a}\overline{\sf a}\overline{\sf a}} &
\frac{1}{{\sf a}\overline{\sf a}} & 0 & 0\smallskip
\\
\frac{{\sf b}{\sf c}}{{\sf a}^4\overline{\sf a}^2}
-
\frac{{\sf e}}{{\sf a}^3\overline{\sf a}}
&
\frac{{\sf c}\overline{\sf b}}{{\sf a}^3
\overline{\sf a}^3}
-
\frac{{\sf d}}{{\sf a}^2\overline{\sf a}^2}
&
\frac{-{\sf c}}{{\sf a}^3\overline{\sf a}^2}
&
\frac{1}{{\sf a}^2\overline{\sf a}}
& 
0
\\
\frac{{\sf b}\overline{\sf c}}{{\sf a}^3\overline{\sf a}^3}
-
\frac{\overline{\sf d}}{{\sf a}^2\overline{\sf a}^2}
&
\frac{\overline{\sf b}\overline{\sf c}}{{\sf a}^2\overline{\sf a}^4}
-
\frac{\overline{\sf e}}{{\sf a}\overline{\sf a}^3}
&
\frac{-\,\overline{\sf c}}{{\sf a}^2\overline{\sf a}^3}
&
0
&
\frac{1}{{\sf a}\overline{\sf a}^2}
\end{array} 
\!\!\right)
\\
&
=:
\left(\!\!
\begin{array}{ccccc}
{\sf a}^\sim & 0 & 0 & 0 & 0
\\
0 & \overline{\sf a}^\sim & 0 & 0 & 0
\\
{\sf b}^\sim & \overline{\sf b}^\sim &
{\sf a}^\sim\overline{\sf a}^\sim & 0 & 0
\\
{\sf e}^\sim & {\sf d}^\sim & {\sf c}^\sim &
{\sf a}^\sim{\sf a}^\sim\overline{\sf a}^\sim & 0
\\
\overline{\sf d}^\sim & \overline{\sf e}^\sim & 
\overline{\sf c}^\sim & 0 & 
{\sf a}^\sim\overline{\sf a}^\sim\overline{\sf a}^\sim
\end{array}
\!\!\right),
\endaligned
\]
also belongs to the subgroup, with:
\[
\aligned
{\sf a}^\sim
&
:=
\frac{1}{\sf a}
\ \ \ \ \ \ \ \ \ \ \ \ \ \ \ \ \ \ \ \ \ \ \ \ \
{\scriptstyle{(\neq\,0)}},
\\
{\sf b}^\sim
&
:=
-\,\frac{\sf b}{{\sf a}{\sf a}\overline{\sf a}},
\\
{\sf c}^\sim
&
:=
-\,\frac{{\sf c}}{{\sf a}{\sf a}{\sf a}\overline{\sf a}\overline{\sf a}},
\\
{\sf d}^\sim
&
:=
\frac{{\sf c}\overline{\sf b}}{
{\sf a}{\sf a}{\sf a}\overline{\sf a}\overline{\sf a}\overline{\sf a}}
-
\frac{{\sf d}}{{\sf a}{\sf a}\overline{\sf a}\overline{\sf a}},
\\
{\sf e}^\sim
&
:=
\frac{{\sf b}{\sf c}}{{\sf a}{\sf a}{\sf a}{\sf a}
\overline{\sf a}\overline{\sf a}}
-
\frac{{\sf e}}{{\sf a}{\sf a}{\sf a}\overline{\sf a}},
\endaligned
\]
which concludes.
\endproof

\noindent{\bf Proposition.}
{\em On a $5$-dimensional CR-generic submanifold:}
\[
\Big(
M^5
\,\subset\,
\C^4
\Big)
\,\,\in\,\,
\text{\sf General Class $\text{\sf III}_{\text{\sf 1}}$},
\]
{\em having biholomorphically invariant $(1, 0)$ CR bundle:}
\[
T^{1,0}M
\,\subset\,
\C\otimes_\R TM,
\]
{\em for any choice of local vector field generator:}
\[
\mathcal{L}
\ \ \ \ \
\text{\rm for}\ \
T^{1,0}M,
\]
{\em the associated frame:}
\[
\aligned
&
\big\{
\mathcal{L},\,\overline{\mathcal{L}},\,
\isqrt\,\big[\mathcal{L},\overline{\mathcal{L}}\big],\,\,
\big[\mathcal{L},\,\isqrt\big[\mathcal{L},
\overline{\mathcal{L}}\big]\big],\,\,
\big[\overline{\mathcal{L}},\,
\isqrt\big[\mathcal{L},\overline{\mathcal{L}}\big]\big]
\big\}
\,=:\,
\\
&
\!\!\!\!\!\!\!\!\!\!\!\!
\,=:\,
\big\{
\mathcal{L},\,\overline{\mathcal{L}},\,\mathcal{T},\,
\mathcal{S},\,\overline{\mathcal{S}}
\big\}
\endaligned
\]
{\em for the full complexified tangent bundle:}
\[
\C\otimes_\R TM
\]
{\em performs a reduction of the full ${\sf GL}_5 ( \C)$-structure
of $\C \otimes_\R TM$ to the $10$-dimensional subgroup:}
\[
{\sf G}_{\text{\sf III}_{\text{\sf 1}}}^{\sf initial}
\,:=\,
\left\{
\left(\!\!
\begin{array}{ccccc}
{\sf a} & 0 & 0 & 0 & 0
\\
0 & \overline{\sf a} & 0 & 0 & 0
\\
{\sf b} & \overline{\sf b} & {\sf a}\overline{\sf a} & 0 & 0
\\
{\sf e} & {\sf d} & {\sf c} & {\sf a}{\sf a}\overline{\sf a} & 0
\\
\overline{\sf d} & \overline{\sf e} & \overline{\sf c} & 0 &
{\sf a}\overline{\sf a}\overline{\sf a}
\end{array}
\!\!\right)
\,\in\,
\mathcal{M}_{5\times 5}(\C)\,\colon\,\,
{\sf a}\,\in\,\C\backslash\{0\},\,\,
{\sf b},\,{\sf c},\,{\sf d},\,{\sf e}\,\in\,\C
\right\}.
\qed
\]


\bigskip

\section{\sf General class $\text{\sf III}_{\text{\sf 2}}$}
\label{general-class-III-2}
\HEAD{\ref{general-class-III-2}.~General class 
$\text{\sf III}_{\text{\sf 2}}$}{
Jo\"el {\sc Merker}, D\'epartement de Math\'ematiques d'Orsay}

\medskip

Equip $\C^4$ with coordinates:
\[
(z,w_1,w_2,w_3)
\,\in\,
\C^4.
\]

Let a connected CR-generic submanifold:
\[
M^5
\subset
\C^4
\]
be of smoothness:
\[
\mathcal{C}^\kappa\ \
(\kappa\geqslant 4),
\ \ \ \ \ \ \ \ \ \ \ \ \ \ \ \ \ \ \ \ \
\text{\rm or}
\ \ \ \ \ \ \ \ \ \ \ \ \ \ \ \ \ \ \ \ \
\mathcal{C}^\infty,
\ \ \ \ \ \ \ \ \ \ \ \ \ \ \ \ \ \ \ \ \
\text{\rm or}
\ \ \ \ \ \ \ \ \ \ \ \ \ \ \ \ \ \ \ \ \
\mathcal{C}^\omega.
\]

Pick a point:
\[
p\in M,
\]
and take a (small) open neighborhood:
\[
p
\in
{\sf U}_p
\subset
\C^4.
\]

Recall (\cite{ Merker-Pocchiola-Sabzevari-5-CR-II,
Merker-5-CR-III}) that if:
\[
\Big(
M^5
\,\subset\,
\C^4
\Big)
\,\,\in\,\,
\text{\sf General Class $\text{\sf III}_{\text{\sf 2}}$},
\]
then:
\[
\aligned
\C\otimes_\R TM
=
T^{1,0}M
+
T^{0,1}M
+
\big[
T^{1,0}M,\,T^{0,1}M\big]
&
+
\big[T^{1,0}M,\,
\big[
T^{1,0}M,\,T^{0,1}M\big]\big]
+
\\
&
+
\big[T^{1,0}M,\,
\big[T^{1,0}M,\,
\big[
T^{1,0}M,\,T^{0,1}M\big]\big]\big].
\endaligned
\]

This means that for any local vector field generator:
\[
\mathcal{L}
\,=\,
\text{\rm section of}\,\,
T^{1,0}\big(M\cap{\sf U}_p\big),
\]
one has at every point $q \in M\cap {\sf U}_p$:
\[
{\bf 5}
=
\rank_\C
\Big(
\mathcal{L}\big\vert_q,\,\,
\overline{\mathcal{L}}\big\vert_q,\,\,
\big[\mathcal{L},\overline{\mathcal{L}}\big]
\Big\vert_q,\,\,
\big[\mathcal{L},\,\big[\mathcal{L},\overline{\mathcal{L}}\big]\big]
\Big\vert_q,\,\,
\big[\mathcal{L},\,\big[\mathcal{L},\,\big[\mathcal{L},
\overline{\mathcal{L}}\big]\big]\big]
\Big\vert_q
\Big).
\]

\medskip

Next, take any (local) biholomorphism:
\[
h\colon\ \ \
{\sf U}_p
\overset{\sim}{\,\longrightarrow\,}
{\sf U}_{p'}'
=
h({\sf U}_p)
\ \ \ \ \ \ \ \ \ \ \ \ \
{\scriptstyle{(p'\,=\,\,h(p))}},
\]
which, when ${\sf U}_p$ is small enough, certainly 
transfers $M \cap {\sf U}_p$ to a certain CR-generic
submanifold:
\[
{M'}^5
\,:=\,
h\big(M\cap{\sf U}_p\big)
\,\subset\,
{\C'}^4,
\]
having of course the same smoothness:
\[
\mathcal{C}^\kappa\ \
(\kappa\geqslant 4),
\ \ \ \ \ \ \ \ \ \ \ \ \ \ \ \ \ \ \ \ \
\text{\rm or}
\ \ \ \ \ \ \ \ \ \ \ \ \ \ \ \ \ \ \ \ \
\mathcal{C}^\infty,
\ \ \ \ \ \ \ \ \ \ \ \ \ \ \ \ \ \ \ \ \
\text{\rm or}
\ \ \ \ \ \ \ \ \ \ \ \ \ \ \ \ \ \ \ \ \
\mathcal{C}^\omega.
\]

\medskip

Take also any local vector field generator:
\[
\mathcal{L}'
\,=\,
\text{\rm section of}\,\,
T^{1,0}M'.
\]
Then necessarily at every point $q' \in M'$, one 
also has (exercise, or {\em see} below):
\[
{\bf 5}
=
\rank_\C
\Big(
\mathcal{L}'\big\vert_{q'},\,\,
\overline{\mathcal{L}}'\big\vert_{q'},\,\,
\big[\mathcal{L}',\overline{\mathcal{L}}'\big]
\Big\vert_{q'},\,\,
\big[\mathcal{L}',\,\big[\mathcal{L}',\overline{\mathcal{L}}'\big]\big]
\Big\vert_{q'},\,\,
\big[\mathcal{L}',\,\big[\mathcal{L}',\,
\big[\mathcal{L}',\overline{\mathcal{L}}'\big]\big]\big]
\Big\vert_{q'}
\Big).
\]

Write the components of $h$ and the target coordinates as:
\[
\!\!\!\!\!\!\!\!\!\!\!\!\!\!\!\!\!\!\!\!
\!\!\!\!\!\!\!\!\!\!\!\!\!\!\!\!\!\!\!\!
\aligned
(z,w_1,w_2,w_3)
&
\,\longmapsto\,
\big(z'(z,w_1,w_2,w_3),\,w_1'(z,w_1,w_2,w_3),\,w_2'(z,w_1,w_2,w_3),\,
w_3'(z,w_1,w_2,w_3)\big)
\\
&\ \ \ \
=:
(z',w_1',w_2',w_3').
\endaligned
\]

Consider the {\em inverse} $h'$
of $h$:
\[
\!\!\!\!\!\!\!\!\!\!\!\!\!\!\!\!\!\!\!\!
\!\!\!\!\!\!\!\!\!\!\!\!\!\!\!\!\!\!\!\!
\aligned
{\sf U}_{p'}'
&
\overset{\sim}{\,\longrightarrow\,}
{\sf U}_p
\ \ \ \ \ \ \ \ \ \ \ \ \ \ \ \ \ \ \ \ \ \ \ \ \ \ \ \ \ \ \ \ \ \ \
\ \ \ \ \ \ \ \ \ \ \ \ \ \ \ \ \ \ \ \ \ \ \ \ \ \ \ \ \ \ \ \ \ \ \
\ \ \ \ \ \ \ \ \ \ \ \ \ \ \ \ \ \ \ \ \ \ \ \ \ \ \ \ \ \ \ \ \ 
{\scriptstyle{(p\,=\,\,h'(p'))}},
\\
(z',w_1',w_2',w_3')
&
\,\longmapsto\,
\big(z(z',w_1',w_2',w_3'),\,w_1(z',w_1',w_2',w_3'),\,
w_2(z',w_1',w_2',w_3'),\,w_3(z',w_1',w_2',w_3')\big)
\\
&\ \ \ \ \,
=
(z,w_1,w_2,w_3).
\endaligned
\]

Then there
exists a nowhere vanishing function:
\[
a\colon\ \ \
M\cap{\sf U}_p
\,\longrightarrow\,
\C\backslash\{0\},
\]
such that:
\[
h_*'
\big(\mathcal{L}'\big)
=
a\,\mathcal{L}.
\]

Simultaneously:
\[
h_*'
\big(\overline{\mathcal{L}}'\big)
=
\overline{a}\,\overline{\mathcal{L}}.
\]

Now setting:
\[
\aligned
\mathcal{T}
&
:=
\isqrt\,
\big[\mathcal{L},\,\overline{\mathcal{L}}\big],
\\
\mathcal{T}'
&
:=
\isqrt\,
\big[\mathcal{L}',\,\overline{\mathcal{L}}'\big],
\endaligned
\]
setting:
\[
\aligned
\mathcal{S}
&
:=
\big[\mathcal{L},\mathcal{T}\big],
\\
\mathcal{S}'
&
:=
\big[\mathcal{L}',\mathcal{T}'\big],
\endaligned
\]
and setting
\[
\aligned
\mathcal{R}
&
=
\big[\mathcal{L},\,
\big[\mathcal{L},\,\mathcal{T}\big]\big],
\\
\mathcal{R}'
&
=
\big[\mathcal{L}',\,
\big[\mathcal{L}',\,\mathcal{T}'\big]\big],
\endaligned
\]
one has two (local) frames:
\[
\aligned
&
\big\{\mathcal{L},\,\overline{\mathcal{L}},\,
\mathcal{T},\,\mathcal{S},\,\mathcal{R}\big\}
\ \ \ \ \ \ \ \ \ \,
\text{\rm for}\ \
\C\otimes_\R T(M\cap{\sf U}_p),
\\
&
\big\{\mathcal{L}',\,\overline{\mathcal{L}}',\,
\mathcal{T}',\,\mathcal{S}',\,\mathcal{R}'\big\}
\ \ \ \ \
\text{\rm for}\ \
\C\otimes_\R TM'.
\endaligned
\]

As for the General Class $\text{\sf III}_{\text{\sf 1}}$:
\[
\aligned
h_*'\big(\mathcal{T}'\big)
&
=
a\overline{a}\,\mathcal{T}
+
\overline{b}\,\overline{\mathcal{L}}
+
b\,\mathcal{L},
\\
h_*'\big(\mathcal{S}'\big)
&
=
aa\overline{a}\,\mathcal{S}
+
c\,\mathcal{T}
+
d\,\overline{\mathcal{L}}
+
e\,\mathcal{L}.
\endaligned
\]

Next, compute:
\[
\aligned
h_*'\big(\mathcal{R}'\big)
&
=
h_*'\big(\big[\mathcal{L}',\mathcal{S}'\big]\big)
\\
&
=
\big[h_*'\big(\mathcal{L}'\big),\,
h_*'\big(\mathcal{S}'\big)\big]
\\
&
=
\big[a\,\mathcal{L},\,\,
aa\overline{a}\,\mathcal{S}
+
c\,\mathcal{T}
+
d\,\overline{\mathcal{L}}
+
e\,\mathcal{L}
\big]
\\
&
=
aaa\overline{a}\,
\underbrace{\big[\mathcal{L},\mathcal{S}\big]}_{
=\,\mathcal{R}}
+
ac\,
\underbrace{\big[\mathcal{L},\mathcal{T}\big]}_{
=\,\mathcal{S}}
+
ad\,
\underbrace{\big[\mathcal{L},\overline{\mathcal{L}}\big]}_{
=\,-\,\isqrt\,\mathcal{T}}
+
ae\,
\zero{\big[\mathcal{L},\mathcal{L}\big]}
+
\\
&
\ \ \ \ \
+
a\mathcal{L}\big(aa\overline{a}\big)
\cdot
\mathcal{S}
+
a\,\mathcal{L}(c)
\cdot
\mathcal{T}
+
a\,\mathcal{L}(d)
\cdot
\overline{\mathcal{L}}
+
a\,\mathcal{L}(e)
\cdot
\mathcal{L}
-
\\
&
\ \ \ \ \
-\,
aa\overline{a}\,\mathcal{S}(a)
\cdot
\mathcal{L}
-
c\,\mathcal{T}(a)
\cdot
\mathcal{L}
-
d\,\overline{\mathcal{L}}(a)
\cdot
\mathcal{L}-
e\,\mathcal{L}(a)
\cdot
\mathcal{L},
\endaligned
\]
that is to say:
\[
\aligned
h_*'\big(\mathcal{R}'\big)
&
=
aaa\overline{a}
\cdot
\mathcal{R}
+
\big(ac+a\,\mathcal{L}(aa\overline{a})\big)
\cdot
\mathcal{S}
+
\\
&
\ \ \ \ \
+
\big(
-\,\isqrt\,ad
+
a\,\mathcal{L}(c)\big)
\cdot
\mathcal{T}
+
\big(a\,\mathcal{L}(d)\big)
\cdot
\overline{\mathcal{L}}
+
\\
&
\ \ \ \ \
+
\big(
a\,\mathcal{L}(e)
-
aa\overline{a}\,\mathcal{S}(a)
-
c\,\mathcal{T}(a)
-
d\,\overline{\mathcal{L}}(a)
-
e\,\mathcal{L}(a)
\big)
\cdot
\mathcal{L}
\\
&
=:
aaa\overline{a}
\cdot
\mathcal{R}
+
f\cdot\mathcal{S}
+
g\cdot\mathcal{T}
+
h\cdot\overline{\mathcal{L}}
+
k\cdot\mathcal{L}.
\endaligned
\]

\medskip\noindent{\bf Summary.}
{\em Through any local biholomorphic equivalences
between CR-generic submanifolds of $\C^4$ belonging
to the General Class} $\text{\sf III}_{\text{\sf 2}}$:
\[
M^5
\overset{\sim}{\,\longrightarrow\,}
{M'}^5,
\]
{\em for any two choices of local vector
field generators:}
\[
\aligned
&
\mathcal{L}
\ \ \ \ \
\text{\rm for}\ \
T^{1,0}M,
\\
&
\mathcal{L}'
\ \ \ \ \
\text{\rm for}\ \
T^{1,0}M',
\endaligned
\]
{\em the transfer of frame obeys the rule:}
\[
\left(\!\!
\begin{array}{c}
\mathcal{L}'
\\
\overline{\mathcal{L}}'
\\
\mathcal{T}'
\\
\mathcal{S}'
\\
\mathcal{R}'
\end{array}
\!\!\right)
=
\left(\!
\begin{array}{ccccc}
a & 0 & 0 & 0 & 0
\\
0 & \overline{a} & 0 & 0 & 0
\\
b & \overline{b} & a\overline{a} & 0 & 0
\\
e & d & c & aa\overline{a} & 0
\\
k & h & g & f & aaa\overline{a}
\end{array}
\!\right)
\left(\!\!
\begin{array}{c}
\mathcal{L}
\\
\overline{\mathcal{L}}
\\
\mathcal{T}
\\
\mathcal{S}
\\
\mathcal{R}
\end{array}
\!\!\right),
\]
{\em for some nine local functions:}
\[
\aligned
a\colon\ \ \
M
&
\,\longrightarrow\,
\C\backslash\{0\},
\\
b,c,d,e,f,g,h,k\colon\ \ \
M
&
\,\longrightarrow\,
\C.
\endaligned
\]

\medskip

This means that the {\em ambiguity} in the choice of a local frame:
\[
\big\{
\mathcal{L},\,
\overline{\mathcal{L}},\,
\mathcal{T},\,
\mathcal{S},\,
\mathcal{R}
\big\}
\ \ \ \ \
\text{\rm for}\ \
\C\otimes_\R TM
\]
which comes\,\,---\,\,{\em naturally from the
point of view of CR geometry}\,\,---\,\,from a choice of a local frame:
\[
\mathcal{L}
\ \ \ \ \
\text{\rm for}\ \
T^{1,0}M
\]
is represented by general changes of frames whose
matrices are of the form: 
\[
\left(\!
\begin{array}{ccccc}
a & 0 & 0 & 0 & 0
\\
0 & \overline{a} & 0 & 0 & 0
\\
b & \overline{b} & a\overline{a} & 0 & 0
\\
e & d & c & aa\overline{a} & 0
\\
k & h & g & f & aaa\overline{a}
\end{array}
\!\right),
\]
the entries being coefficient-functions depending on existing
equivalences $M \overset{ \sim}{ \longrightarrow} M'$,
and\big/or on the choice of local coordinates.

\medskip

Thus, the {\em initial
$G$-structure} for the biholomorphic equivalence problem between
CR-generic submanifolds $M^5 \subset \C^4$ belonging to the General Class
$\text{\sf III}_{\text{\sf 2}}$ is a reduction of the full linear group 
${\sf GL}_5 ( \C)$ to
the mentioned subgroup in which the\,\,---\,\,possibly
unknown\,\,---\,\,functions $a$, $b$,
$c$, $d$, $e$ are {\em replaced} by independent
complex variables.

\medskip\noindent{\bf Lemma.}
{\em The set of matrices:}
\[
\!\!\!\!\!\!\!\!\!\!\!\!\!\!\!\!\!\!\!\!
\!\!\!\!\!\!\!\!\!\!\!\!\!\!\!
{\sf G}_{\text{\sf III}_{\text{\sf 2}}}^{\sf initial}
\,:=\,
\left\{
\left(\!\!
\begin{array}{ccccc}
{\sf a} & 0 & 0 & 0 & 0
\\
0 & \overline{\sf a} & 0 & 0 & 0
\\
{\sf b} & \overline{\sf b} & {\sf a}\overline{\sf a} & 0 & 0
\\
{\sf e} & {\sf d} & {\sf c} & {\sf a}{\sf a}\overline{\sf a} & 0
\\
{\sf k} & {\sf h} & {\sf g} & {\sf f} &
{\sf a}{\sf a}{\sf a}\overline{\sf a}
\end{array}
\!\!\right)
\,\in\,
\mathcal{M}_{5\times 5}(\C)\,\colon\,\,
{\sf a}\,\in\,\C\backslash\{0\},\,\,
{\sf b},\,{\sf c},\,{\sf d},\,{\sf e}\,
{\sf f},\,{\sf g},\,{\sf h},\,{\sf k}\in\,\C
\right\}
\]
{\em is a (closed) $18$-dimensional
real matrix subgroup of the full:}
\[
{\sf GL}_5(\C)
\,:=\,
\left\{
\pi
=
\left(\!\!
\begin{array}{ccccc}
\pi_{1,1} & \pi_{1,2} & \pi_{1,3} & \pi_{1,4} & \pi_{1,5}
\\
\pi_{2,1} & \pi_{2,2} & \pi_{3,3} & \pi_{2,4} & \pi_{2,5}
\\
\pi_{3,1} & \pi_{3,2} & \pi_{3,3} & \pi_{3,4} & \pi_{3,5}
\\
\pi_{4,1} & \pi_{4,2} & \pi_{4,3} & \pi_{4,4} & \pi_{4,5}
\\
\pi_{5,1} & \pi_{5,2} & \pi_{5,3} & \pi_{5,4} & \pi_{5,5}
\end{array}
\!\!\right)
\,\in\,
\mathcal{M}_{5\times 5}(\C)\,\colon\,\,
0
\neq
\det\,\pi
\right\}.
\]

\proof
Closedness under multiplication (composition):
\[
\aligned
&
\left(\!\!
\begin{array}{ccccc}
{\sf a}_1 & 0 & 0 & 0 & 0
\\
0 & \overline{\sf a}_1 & 0 & 0 & 0
\\
{\sf b}_1 & \overline{\sf b}_1 &
{\sf a}_1\overline{\sf a}_1 & 0 & 0
\\
{\sf e}_1 & {\sf d}_1 & {\sf c}_1 & {\sf a}_1{\sf a}_1
\overline{\sf a}_1 & 0
\\
{\sf k}_1 & {\sf h}_1 & {\sf g}_1 & {\sf f}_1 &
{\sf a}_1{\sf a}_1{\sf a}_1\overline{\sf a}_1
\end{array}
\!\!\right)
\cdot
\left(\!\!
\begin{array}{ccccc}
{\sf a}_2 & 0 & 0 & 0 & 0
\\
0 & \overline{\sf a}_2 & 0 & 0 & 0
\\
{\sf b}_2 & \overline{\sf b}_2 &
{\sf a}_2\overline{\sf a}_2 & 0 & 0
\\
{\sf e}_2 & {\sf d}_2 & {\sf c}_2 & {\sf a}_2{\sf a}_2
\overline{\sf a}_2 & 0
\\
{\sf k}_2 & {\sf h}_2 & {\sf g}_2 & {\sf f}_2 &
{\sf a}_2{\sf a}_2{\sf a}_2\overline{\sf a}_2
\end{array}
\!\!\right)
\\
&
\ \ \ \ \ \ \ \ \ \
=
\left(\!\!
\begin{array}{ccccc}
{\sf a}_1{\sf a}_2 & 0 & 0 & 0 & 0
\\
0 & \overline{\sf a}_1\overline{\sf a}_2 & 0 & 0 & 0
\\
{\sf b}_1{\sf a}_2
\!+\!{\sf a}_1\overline{\sf a}_1{\sf b}_2 & 
\overline{\sf b}_1\overline{\sf a}_2
\!+\!{\sf a}_1\overline{\sf a}_1\overline{\sf b}_2 &
{\sf a}_1\overline{\sf a}_1{\sf a}_2\overline{\sf a}_2 & 0 & 0\medskip
\\
\substack{{\sf e}_1{\sf a}_2+{\sf c}_1{\sf b}_2+\\
+{\sf a}_1{\sf a}_1\overline{\sf a}_1{\sf e}_2}
&
\substack{{\sf d}_1\overline{\sf a}_2+{\sf c}_1\overline{\sf b}_2+\\
+{\sf a}_1{\sf a}_1\overline{\sf a}_1{\sf d}_2}
&
\substack{{\sf c}_1{\sf a}_2\overline{\sf a}_2+\\
{\sf a}_1{\sf a}_1\overline{\sf a}_1{\sf c}_2}
&
\substack{
{\sf a}_1{\sf a}_1\overline{\sf a}_1
{\sf a}_2{\sf a}_2\overline{\sf a}_2} & 0\medskip
\\
\substack{{\sf k}_1{\sf a}_2+{\sf g}_1{\sf b}_2+\\
+{\sf f}_1{\sf e}_2
+{\sf a}_1{\sf a}_1{\sf a}_1\overline{\sf a}_1{\sf k}_2}
&
\substack{
{\sf h}_1\overline{\sf a}_2+{\sf g}_1\overline{\sf b}_2+\\
+{\sf f}_1{\sf d}_2
+{\sf a}_1{\sf a}_1{\sf a}_1\overline{\sf a}_1{\sf h}_2}
&
\substack{
{\sf g}_1{\sf a}_2\overline{\sf a}_2+{\sf f}_1{\sf c}_2\\
+{\sf a}_1{\sf a}_1{\sf a}_1\overline{\sf a}_1{\sf g}_2}
&
\substack{
{\sf f}_1{\sf a}_2{\sf a}_2\overline{\sf a}_2+\\
+{\sf a}_1{\sf a}_1{\sf a}_1\overline{\sf a}_1{\sf f}_2}
&
\substack{
{\sf a}_1{\sf a}_1\overline{\sf a}_1\overline{\sf a}_1\cdot\\
\cdot{\sf a}_2{\sf a}_2\overline{\sf a}_2\overline{\sf a}_2}
\end{array}
\!\!\right)
\\
&
\ \ \ \ \ \ \ \ \ \
=:
\left(\!\!
\begin{array}{ccccc}
{\sf a}_3 & 0 & 0 & 0 & 0
\\
0 & \overline{\sf a}_3 & 0 & 0 & 0
\\
{\sf b}_3 & \overline{\sf b}_3 &
{\sf a}_3\overline{\sf a}_3 & 0 & 0
\\
{\sf e}_3 & {\sf d}_3 & {\sf c}_3 & {\sf a}_3{\sf a}_3\overline{\sf a}_3 & 0
\\
{\sf k}_3 & {\sf h}_3 & {\sf g}_3 & {\sf f}_3 &
{\sf a}_3{\sf a}_3{\sf a}_3\overline{\sf a}_3
\end{array}
\!\!\right)
\endaligned
\]
is visibly clear, after setting:
\[
\aligned
{\sf a}_3
&
:=
{\sf a}_1{\sf a}_2
\ \ \ \ \ \ \ \ \ \ \ \ \ \ \ \ \ \ \ \ \ \ \ \ \
{\scriptstyle{(\neq\,0,\,\,
\text{\rm again})}},
\\
{\sf b}_3
&
:=
{\sf b}_1{\sf a}_2
+
{\sf a}_1\overline{\sf a}_1{\sf b}_2,
\\
{\sf c}_3
&
:=
{\sf c}_1{\sf a}_2\overline{\sf a}_2
+
{\sf a}_1{\sf a}_1\overline{\sf a}_1{\sf c}_2,
\\
{\sf d}_3
&
:=
{\sf d}_1\overline{\sf a}_2
+
{\sf c}_1\overline{\sf b}_2
+
{\sf a}_1{\sf a}_1\overline{\sf a}_1{\sf d}_2,
\\
{\sf e}_3
&
:=
{\sf e}_1{\sf a}_2
+
{\sf c}_1{\sf b}_2
+
{\sf a}_1{\sf a}_1\overline{\sf a}_1{\sf e}_2,
\\
{\sf f}_3
&
:=
{\sf f}_1{\sf a}_2{\sf a}_2\overline{\sf a}_2
+
{\sf a}_1{\sf a}_1{\sf a}_1\overline{\sf a}_1{\sf f}_2,
\\
{\sf g}_3
&
:=
{\sf g}_1{\sf a}_2\overline{\sf a}_2
+
{\sf f}_1{\sf c}_2
+
{\sf a}_1{\sf a}_1{\sf a}_1\overline{\sf a}_1{\sf g}_2,
\\
{\sf h}_3
&
:=
{\sf h}_1\overline{\sf a}_2
+
{\sf g}_1\overline{\sf b}_2
+
{\sf f}_1{\sf d}_2
+
{\sf a}_1{\sf a}_1{\sf a}_1\overline{\sf a}_1{\sf h}_2,
\\
{\sf k}_3
&
:=
{\sf k}_1{\sf a}_2
+
{\sf g}_1{\sf b}_2
+
{\sf f}_1{\sf e}_2
+
{\sf a}_1{\sf a}_1{\sf a}_1\overline{\sf a}_1{\sf k}_2.
\endaligned
\]

Quite similarly, the inverse:
\[
\aligned
\left(\!\!
\begin{array}{ccccc}
{\sf a} & 0 & 0 & 0 & 0
\\
0 & \overline{\sf a} & 0 & 0 & 0
\\
{\sf b} & \overline{\sf b} &
{\sf a}\overline{\sf a} & 0 & 0
\\
{\sf e} & {\sf d} & {\sf c} & 
{\sf a}{\sf a}\overline{\sf a} & 0
\\
{\sf k} & {\sf h} & {\sf g} & {\sf f} &
{\sf a}{\sf a}{\sf a}\overline{\sf a}
\end{array}
\!\!\right)^{-1}
&
=
\left(\!\!
\begin{array}{ccccc}
\frac{1}{\sf a} & 0 & 0 & 0 & 0
\\
0 & \frac{1}{\overline{\sf a}} & 0 & 0 & 0\smallskip
\\
\frac{-{\sf b}}{{\sf a}{\sf a}\overline{\sf a}} & 
\frac{-\overline{\sf b}}{{\sf a}\overline{\sf a}\overline{\sf a}} &
\frac{1}{{\sf a}\overline{\sf a}} & 0 & 0\smallskip
\\
\frac{{\sf b}{\sf c}}{{\sf a}^4\overline{\sf a}^2}
-
\frac{{\sf e}}{{\sf a}^3\overline{\sf a}}
&
\frac{{\sf c}\overline{\sf b}}{{\sf a}^3
\overline{\sf a}^3}
-
\frac{{\sf d}}{{\sf a}^2\overline{\sf a}^2}
&
\frac{-{\sf c}}{{\sf a}^3\overline{\sf a}^2}
&
\frac{1}{{\sf a}^2\overline{\sf a}}
& 
0\bigskip
\\
\substack{
\frac{-{\sf b}{\sf c}{\sf f}}{{\sf a}^7\overline{\sf a}^3}
+
\frac{{\sf b}{\sf g}}{{\sf a}^5\overline{\sf a}^2}+\\
+
\frac{{\sf e}{\sf f}}{{\sf a}^6\overline{\sf a}^2}
-
\frac{{\sf k}}{{\sf a}^4\overline{\sf a}}}
&
\substack{
\frac{-{\sf f}{\sf c}\overline{\sf b}}{{\sf a}^6\overline{\sf a}^4}
+
\frac{{\sf b}\overline{\sf b}}{{\sf a}^4\overline{\sf a}^3}+\\
+\frac{{\sf f}{\sf d}}{{\sf a}^5\overline{\sf a}^3}
-
\frac{{\sf h}}{{\sf a}^3\overline{\sf a}^2}}
&
\frac{{\sf f}{\sf c}}{{\sf a}^6\overline{\sf a}^3}
-
\frac{{\sf g}}{{\sf a}^4\overline{\sf a}^2}
&
\frac{-\,{\sf f}}{{\sf a}^5\overline{\sf a}^2}
&
\frac{1}{{\sf a}^3\overline{\sf a}}
\end{array} 
\!\!\right)
\\
&
=:
\left(\!\!
\begin{array}{ccccc}
{\sf a}^\sim & 0 & 0 & 0 & 0
\\
0 & \overline{\sf a}^\sim & 0 & 0 & 0
\\
{\sf b}^\sim & \overline{\sf b}^\sim &
{\sf a}^\sim\overline{\sf a}^\sim & 0 & 0
\\
{\sf e}^\sim & {\sf d}^\sim & {\sf c}^\sim &
{\sf a}^\sim{\sf a}^\sim\overline{\sf a}^\sim & 0
\\
{\sf k}^\sim & {\sf h}^\sim & {\sf g}^\sim & {\sf f}^\sim &
{\sf a}^\sim{\sf a}^\sim{\sf a}^\sim\overline{\sf a}^\sim
\end{array}
\!\!\right),
\endaligned
\]
also belongs to the subgroup, with:
\[
\aligned
{\sf a}^\sim
&
:=
\frac{1}{\sf a}
\ \ \ \ \ \ \ \ \ \ \ \ \ \ \ \ \ \ \ \ \ \ \ \ \
{\scriptstyle{(\neq\,0)}},
\\
{\sf b}^\sim
&
:=
-\,\frac{\sf b}{{\sf a}{\sf a}\overline{\sf a}},
\\
{\sf c}^\sim
&
:=
-\,\frac{{\sf c}}{{\sf a}{\sf a}{\sf a}\overline{\sf a}\overline{\sf a}},
\\
{\sf d}^\sim
&
:=
\frac{{\sf c}\overline{\sf b}}{
{\sf a}{\sf a}{\sf a}\overline{\sf a}\overline{\sf a}\overline{\sf a}}
-
\frac{{\sf d}}{{\sf a}{\sf a}\overline{\sf a}\overline{\sf a}},
\\
{\sf e}^\sim
&
:=
\frac{{\sf b}{\sf c}}{{\sf a}{\sf a}{\sf a}{\sf a}
\overline{\sf a}\overline{\sf a}}
-
\frac{{\sf e}}{{\sf a}{\sf a}{\sf a}\overline{\sf a}},
\\
{\sf f}^\sim
&
:=
-\,\frac{{\sf f}}{{\sf a}^5\overline{\sf a}^2},
\\
{\sf g}^\sim
&
:=
\frac{{\sf c}{\sf f}}{{\sf a}^6\overline{\sf a}^3}
-
\frac{{\sf g}}{{\sf a}^4\overline{\sf a}^2},
\\
{\sf h}^\sim
&
:=
-\,
\frac{{\sf f}{\sf c}\overline{\sf b}}{{\sf a}^6\overline{\sf a}^4}
+
\frac{{\sf g}\overline{\sf b}}{{\sf a}^4\overline{\sf a}^3}
+
\frac{{\sf f}{\sf d}}{{\sf a}^5\overline{\sf a}^3}
-
\frac{{\sf h}}{{\sf a}^3\overline{\sf a}^2},
\\
{\sf k}^\sim
&
:=
-\,
\frac{{\sf b}{\sf c}{\sf f}}{{\sf a}^7\overline{\sf a}^3}
+
\frac{{\sf b}{\sf g}}{{\sf a}^5\overline{\sf a}^2}
+
\frac{{\sf e}{\sf f}}{{\sf a}^6\overline{\sf a}^2}
-
\frac{{\sf k}}{{\sf a}^4\overline{\sf a}^1},
\endaligned
\]
which concludes.
\endproof

\noindent{\bf Proposition.}
{\em On a $5$-dimensional CR-generic submanifold:}
\[
\Big(
M^5
\,\subset\,
\C^4
\Big)
\,\,\in\,\,
\text{\sf General Class $\text{\sf III}_{\text{\sf 2}}$},
\]
{\em having biholomorphically invariant $(1, 0)$ CR bundle:}
\[
T^{1,0}M
\,\subset\,
\C\otimes_\R TM,
\]
{\em for any choice of local vector field generator:}
\[
\mathcal{L}
\ \ \ \ \
\text{\rm for}\ \
T^{1,0}M,
\]
{\em the associated frame:}
\[
\aligned
&
\big\{
\mathcal{L},\,\overline{\mathcal{L}},\,
\isqrt\,\big[\mathcal{L},\overline{\mathcal{L}}\big],\,\,
\big[\mathcal{L},\,\isqrt\big[\mathcal{L},
\overline{\mathcal{L}}\big]\big],\,\,
\big[\mathcal{L},\,\big[\mathcal{L},\,\isqrt\big[\mathcal{L},
\overline{\mathcal{L}}\big]\big]\big]
\big\}
\,=:\,
\\
&
\!\!\!\!\!\!\!\!\!\!\!\!
\,=:\,
\big\{
\mathcal{L},\,\overline{\mathcal{L}},\,\mathcal{T},\,
\mathcal{S},\,\mathcal{R}
\big\}
\endaligned
\]
{\em for the full complexified tangent bundle:}
\[
\C\otimes_\R TM
\]
{\em performs a reduction of the full ${\sf GL}_5 ( \C)$-structure
of $\C \otimes_\R TM$ to the $18$-dimensional subgroup:}
\[
\!\!\!\!\!\!\!\!\!\!\!\!\!\!\!\!\!\!\!\!
\!\!\!\!\!\!\!\!\!\!\!\!\!\!\!\!\!\!\!\!
{\sf G}_{\text{\sf III}_{\text{\sf 2}}}^{\sf initial}
\,:=\,
\left\{
\left(\!\!
\begin{array}{ccccc}
{\sf a} & 0 & 0 & 0 & 0
\\
0 & \overline{\sf a} & 0 & 0 & 0
\\
{\sf b} & \overline{\sf b} & {\sf a}\overline{\sf a} & 0 & 0
\\
{\sf e} & {\sf d} & {\sf c} & {\sf a}{\sf a}\overline{\sf a} & 0
\\
{\sf k} & {\sf h} & {\sf g} & {\sf f} &
{\sf a}{\sf a}{\sf a}\overline{\sf a}
\end{array}
\!\!\right)
\,\in\,
\mathcal{M}_{5\times 5}(\C)\,\colon\,\,
{\sf a}\,\in\,\C\backslash\{0\},\,\,
{\sf b},\,{\sf c},\,{\sf d},\,{\sf e},\,
{\sf f},\,{\sf g},\,{\sf h},\,{\sf k}\,
\,\in\,\C
\right\}.
\qed
\]


\bigskip

\section{\sf General class $\text{\sf IV}_{\text{\sf 1}}$}
\label{general-class-IV-1}
\HEAD{\ref{general-class-IV-1}.~General class 
$\text{\sf IV}_{\text{\sf 1}}$}{
Jo\"el {\sc Merker}, D\'epartement de Math\'ematiques d'Orsay}

\medskip

Equip $\C^3$ with coordinates:
\[
(z_1,z_2,w)
\,\in\,
\C^3.
\]

Let a connected hypersurface:
\[
M^5
\,\subset\,
\C^3
\]
be of smoothness:

\[
\mathcal{C}^\kappa\ \
(\kappa\geqslant 3),
\ \ \ \ \ \ \ \ \ \ \ \ \ \ \ \ \ \ \ \ \
\text{\rm or}
\ \ \ \ \ \ \ \ \ \ \ \ \ \ \ \ \ \ \ \ \
\mathcal{C}^\infty,
\ \ \ \ \ \ \ \ \ \ \ \ \ \ \ \ \ \ \ \ \
\text{\rm or}
\ \ \ \ \ \ \ \ \ \ \ \ \ \ \ \ \ \ \ \ \
\mathcal{C}^\omega.
\]

Pick a point:
\[
p\in M,
\]
and take a (small) open neighborhood:
\[
p
\in
{\sf U}_p
\subset
\C^3.
\]
Take a local frame:
\[
\big\{
\mathcal{L}_1,\,\mathcal{L}_2
\big\}
\ \ \ \ \
\text{\rm for}\ \
T^{1,0}\big(M\cap{\sf U}_p\big).
\]

Recall (\cite{ Merker-Pocchiola-Sabzevari-5-CR-II,
Merker-5-CR-III}) that:
\[
\Big(
M^5
\,\subset\,
\C^4
\Big)
\,\,\in\,\,
\text{\sf General Class $\text{\sf IV}_{\text{\sf 1}}$},
\]
if firstly, after a possible constant change of frame:
\[
\left(\!\!
\begin{array}{c}
\mathcal{L}_1
\\
\mathcal{L}_2
\end{array}
\!\!\right)
\,\longmapsto\,
\underbrace{\left(\!\!
\begin{array}{cc}
\alpha & \beta
\\
\gamma & \delta
\end{array}
\!\!\right)}_{\in\,{\sf GL}_2(\C)}
\left(\!\!
\begin{array}{c}
\mathcal{L}_1
\\
\mathcal{L}_2
\end{array}
\!\!\right),
\]
setting:
\[
\mathcal{T}
\,:=\,
\isqrt\,\big[\mathcal{L}_1,
\overline{\mathcal{L}}_1\big],
\]
the five field:
\[
\big\{\mathcal{L}_1,\,\mathcal{L}_2,\,
\overline{\mathcal{L}}_1,\,\overline{\mathcal{L}}_2,\,
\mathcal{T}
\big\}
\]
make up a (local) frame for:
\[
\C\otimes_\R TM,
\]
and if, secondly, the Levi form of $M$ is nondegenerate
(of rank $2$) at every point, a second condition 
which reads as follows.

\smallskip

Alltogether:
\[
\aligned
\isqrt\,\big[\mathcal{L}_1,\,\overline{\mathcal{L}}_1\big]
&
\,=\,
\mathcal{T},
\\
\isqrt\,\big[\mathcal{L}_2,\,\overline{\mathcal{L}}_1\big]
&
\,=\,
A\cdot\mathcal{T}
+
B_1\cdot\overline{\mathcal{L}}_1
+
B_2\cdot\overline{\mathcal{L}}_2
+
D_1\cdot\mathcal{L}_1
+
D_2\cdot\mathcal{L}_2,
\\
\isqrt\,\big[\mathcal{L}_1,\,\overline{\mathcal{L}}_2\big]
&
\,=\,
\overline{A}\cdot\mathcal{T}
+
\overline{D}_1\cdot\overline{\mathcal{L}}_1
+
\overline{D}_2\cdot\overline{\mathcal{L}}_2
+
\overline{B}_1\cdot\mathcal{L}_1
+
\overline{B}_2\cdot\mathcal{L}_2,
\\
\isqrt\,\big[\mathcal{L}_2,\,\overline{\mathcal{L}}_2\big]
&
\,=\,
C\cdot\mathcal{T}
+
\overline{E}_1\cdot\overline{\mathcal{L}}_1
+
\overline{E}_2\cdot\overline{\mathcal{L}}_2
+
E_1\cdot\mathcal{L}_1
+
E_2\cdot\mathcal{L}_2,
\endaligned
\]
for certain functions:
\[
\aligned
C\,\,
\colon\ \ \
M\cap{\sf U}_p
&
\,\longrightarrow\,
\R,
\\
A,\,B_1,\,B_2,\,D_1,\,D_2,\,E_1,\,E_2\,\,
\colon\ \ \
M\cap{\sf U}_p
&
\,\longrightarrow\,
\C.
\endaligned
\]

Then the Levi form of $M$ is of rank $2$ at every point
$q \in M \cap {\sf U}_p$ if and only if
(mental exercise, {\em cf.}~\cite{ Merker-Pocchiola-Sabzevari-5-CR-II}):
\[
0
\,\neq\,
\det\,
\left(\!
\begin{array}{cc}
1 & A 
\\
\overline{A} & C
\end{array}
\!\right)(q)
\ \ \ \ \ \ \ \ \ \ \ \ \
{\scriptstyle{(\forall\,q\,\in\,M\cap{\sf U}_p)}}.
\]

Next, take any (local) biholomorphism:
\[
h\colon\ \ \
{\sf U}_p
\overset{\sim}{\,\longrightarrow\,}
{\sf U}_{p'}'
=
h({\sf U}_p)
\ \ \ \ \ \ \ \ \ \ \ \ \
{\scriptstyle{(p'\,=\,\,h(p))}},
\]
which, when ${\sf U}_p$ is small enough, certainly 
transfers $M \cap {\sf U}_p$ to a certain hypersurface:
\[
{M'}^5
\,:=\,
h\big(M\cap{\sf U}_p\big)
\,\subset\,
{\C'}^3,
\]
having the same smoothness (mental exercise):
\[
\mathcal{C}^\kappa\ \
(\kappa\geqslant 3),
\ \ \ \ \ \ \ \ \ \ \ \ \ \ \ \ \ \ \ \ \
\text{\rm or}
\ \ \ \ \ \ \ \ \ \ \ \ \ \ \ \ \ \ \ \ \
\mathcal{C}^\infty,
\ \ \ \ \ \ \ \ \ \ \ \ \ \ \ \ \ \ \ \ \
\text{\rm or}
\ \ \ \ \ \ \ \ \ \ \ \ \ \ \ \ \ \ \ \ \
\mathcal{C}^\omega.
\]
Consider the inverse of $h$:
\[
h'\colon\ \ \
{\sf U}_{p'}'
\overset{\sim}{\,\longrightarrow\,}
{\sf U}_p
\ \ \ \ \ \ \ \ \ \ \ \ \ \ \ \ \ \ \ \ \ \ \ \ \ \ \ \ \ \
{\scriptstyle{(p\,=\,\,h'(p'))}}.
\]

\medskip

Take also any two local vector field generators:
\[
\big\{
\mathcal{L}_1',\,
\mathcal{L}_2'
\big\}
\]
making a frame for $T^{1, 0} M'$. After a possible
constant ${\sf GL}_2 ( \C)$-multiplication, setting similarly:
\[
\mathcal{T}'
:=
\isqrt\,\big[\mathcal{L}_1',\overline{\mathcal{L}}_1'\big],
\]
the five fields:
\[
\big\{
\mathcal{L}_1',\,\mathcal{L}_2',\,
\overline{\mathcal{L}}_1',\,\overline{\mathcal{L}}_2',\,
\mathcal{T}'
\big\}
\]
also make up (local) frame for:
\[
\C\otimes_\R TM'.
\]

Then (\cite{ Merker-Pocchiola-Sabzevari-5-CR-II}), 
there are certain multiplicator functions with:
\[
\aligned
h_*'\big(\mathcal{L}_1'\big)
&
=
a_{11}\,\mathcal{L}_1
+
a_{21}\,\mathcal{L}_2,
\\
h_*'\big(\mathcal{L}_2'\big)
&
=
a_{12}\,\mathcal{L}_1
+
a_{22}\,\mathcal{L}_2,
\endaligned
\]
with:
\[
0
\,\neq\,
\det\,
\left(\!
\begin{array}{cc}
a_{11} & a_{12} 
\\
a_{21} & a_{22}
\end{array}
\!\right)(q)
\ \ \ \ \ \ \ \ \ \ \ \ \
{\scriptstyle{(\forall\,q\,\in\,M\cap{\sf U}_p)}},
\]
whence also:
\[
\aligned
h_*'\big(\overline{\mathcal{L}}_1'\big)
&
=
\overline{a}_{11}\,\overline{\mathcal{L}}_1
+
\overline{a}_{21}\,\overline{\mathcal{L}}_2,
\\
h_*'\big(\overline{\mathcal{L}}_2'\big)
&
=
\overline{a}_{12}\,\overline{\mathcal{L}}_1
+
\overline{a}_{22}\,\overline{\mathcal{L}}_2.
\endaligned
\]

\smallskip

Now, compute:
\[
\aligned
h_*'\big(\mathcal{T}'\big)
&
=
h_*'\big(\isqrt\,\big[\mathcal{L}_1',\overline{\mathcal{L}}_1'\big]\big)
\\
&
=
\isqrt\,
\big[
h_*'\big(\mathcal{L}_1'\big),\,
h_*'\big(\overline{\mathcal{L}}_1'\big)
\big]
\\
&
=
\big[a_{11}\,\mathcal{L}_1+a_{21}\,\mathcal{L}_2,\,\,
\overline{a}_{11}\,\overline{\mathcal{L}}_1
+
\overline{a}_{21}\,\overline{\mathcal{L}}_2
\big]
\\
&
=
a_{11}\overline{a}_{11}\,\isqrt\,
\underbrace{\big[\mathcal{L}_1,\overline{\mathcal{L}}_1\big]}_{
=\,\mathcal{T}}
+
a_{21}\overline{a}_{11}\,\isqrt\,
\big[\mathcal{L}_2,\overline{\mathcal{L}}_1\big]
+
\\
&
\ \ \ \ \
+
a_{11}\overline{a}_{21}\,\isqrt\,
\big[\mathcal{L}_1,\overline{\mathcal{L}}_2\big]
+
a_{21}\overline{a}_{21}\,\isqrt\,
\big[\mathcal{L}_2,\overline{\mathcal{L}}_2\big]
+
\\
&
\ \ \ \ \
+
\isqrt\,\big(
a_{11}\,\mathcal{L}_1\big(\overline{a}_{11}\big)
+
a_{21}\,\mathcal{L}_2\big(\overline{a}_{11}\big)
\big)
\cdot
\overline{\mathcal{L}}_1
+
\\
&
\ \ \ \ \
+
\isqrt\,\big(
a_{11}\,\mathcal{L}_1\big(\overline{a}_{21}\big)
+
a_{21}\,\mathcal{L}_2\big(\overline{a}_{21}\big)
\cdot
\overline{\mathcal{L}}_2
-
\\
&
\ \ \ \ \
-
\isqrt\,\big(
\overline{a}_{11}\,\overline{\mathcal{L}}_1(a_{11})
+
\overline{a}_{21}\,\overline{\mathcal{L}}_2(a_{11})
\big)
\cdot
\mathcal{L}_1
-
\\
&
\ \ \ \ \
-
\isqrt\,\big(
\overline{a}_{11}\,\overline{\mathcal{L}}_1(a_{21})
+
\overline{a}_{21}\,\overline{\mathcal{L}}_2(a_{21})
\big)
\cdot
\mathcal{L}_2,
\endaligned
\]
which, taking account of what precedes, is:
\[
\!\!\!\!\!\!\!\!\!\!\!\!\!\!\!\!\!\!\!\!
\aligned
h_*'\big(\mathcal{T}'\big)
&
=
\Big(
a_{11}\overline{a}_{11}
+
a_{21}\overline{a}_{11}A
+
a_{11}\overline{a}_{21}\overline{A}
+
a_{21}\overline{a}_{21}C
\Big)
\cdot
\mathcal{T}
+
\\
&
\ \ \ \ \
+
\Big(
a_{21}\overline{a}_{11}B_1
+
a_{11}\overline{a}_{21}\overline{D}_1
+
a_{21}\overline{a}_{21}\overline{E}_1
+
\isqrt\,a_{11}\mathcal{L}_1\big(\overline{a}_{11}\big)
+
\isqrt\,a_{21}\mathcal{L}_2\big(\overline{a}_{11}\big)
\Big)
\cdot
\overline{\mathcal{L}}_1
+
\\
&
\ \ \ \ \
+
\Big(
a_{21}\overline{a}_{11}B_1
+
a_{11}\overline{a}_{21}\overline{D}_2
+
a_{21}\overline{a}_{21}\overline{E}_2
+
\isqrt\,a_{11}\mathcal{L}_1\big(\overline{a}_{21}\big)
+
\isqrt\,a_{21}\mathcal{L}_2\big(\overline{a}_{21}\big)
\Big)
\cdot
\overline{\mathcal{L}}_2
+
\\
&
\ \ \ \ \
+
\Big(
a_{21}\overline{a}_{11}D_1
+
a_{11}\overline{a}_{21}\overline{B}_1
+
a_{21}\overline{a}_{21}E_1
-
\isqrt\,\overline{a}_{11}\overline{\mathcal{L}}_1(a_{11})
-
\isqrt\,\overline{a}_{21}\overline{\mathcal{L}}_2(a_{11})
\Big)
\cdot
\mathcal{L}_1
+
\\
&
\ \ \ \ \
+
\Big(
a_{21}\overline{a}_{11}D_2
+
a_{11}\overline{a}_{21}\overline{B}_2
+
a_{21}\overline{a}_{21}E_2
-
\isqrt\,\overline{a}_{11}\overline{\mathcal{L}}_1(a_{21})
-
\isqrt\,\overline{a}_{21}\overline{\mathcal{L}}_2(a_{21})
\Big)
\cdot
\mathcal{L}_2.
\endaligned
\]

The first appearing coefficient-function:
\[
c
:=
a_{11}\overline{a}_{11}
+
a_{21}\overline{a}_{11}A
+
a_{11}\overline{a}_{21}\overline{A}
+
a_{21}\overline{a}_{21}C
\]
is real-valued and vanishes nowhere (mental exercise):
\[
c\colon\ \ \
M\cap{\sf U}_p
\,\longrightarrow\,
\R\backslash\{0\}.
\]
The two functions (and their conjugates):
\[
\aligned
b_1
&
:=
a_{21}\overline{a}_{11}D_1
+
a_{11}\overline{a}_{21}\overline{B}_1
+
a_{21}\overline{a}_{21}E_1
-
\isqrt\,\overline{a}_{11}\overline{\mathcal{L}}_1(a_{11})
-
\isqrt\,\overline{a}_{21}\overline{\mathcal{L}}_2(a_{11}),
\\
b_2
&
:=
a_{21}\overline{a}_{11}D_2
+
a_{11}\overline{a}_{21}\overline{B}_2
+
a_{21}\overline{a}_{21}E_2
-
\isqrt\,\overline{a}_{11}\overline{\mathcal{L}}_1(a_{21})
-
\isqrt\,\overline{a}_{21}\overline{\mathcal{L}}_2(a_{21}),
\endaligned
\]
are complex-valued:
\[
\aligned
b_1,\,\,b_2\,\,\colon\ \ \ 
M\cap{\sf U}_p
\,\longrightarrow\,
\C.
\endaligned
\]

\medskip\noindent{\bf Summary.}
{\em Through any local biholomorphic equivalences
between CR-generic submanifolds of $\C^3$ belonging
to the General Class} $\text{\sf IV}_{\text{\sf 1}}$:
\[
M^5
\overset{\sim}{\,\longrightarrow\,}
{M'}^5,
\]
{\em for any two choices of pairs of local vector
field generators:}
\[
\aligned
&
\big\{\mathcal{L}_1,\mathcal{L}_2\big\}
\ \ \ \ \
\text{\rm for}\ \
T^{1,0}M,
\\
&
\big\{\mathcal{L}_1',\mathcal{L}_2'\big\}
\ \ \ \ \
\text{\rm for}\ \
T^{1,0}M',
\endaligned
\]
{\em the transfer of frame obeys the rule:}
\[
\left(\!\!
\begin{array}{c}
\mathcal{L}_1'
\\
\mathcal{L}_2'
\\
\overline{\mathcal{L}}_1'
\\
\overline{\mathcal{L}}_2'
\\
\mathcal{T}'
\end{array}
\!\!\right)
=
\left(\!
\begin{array}{ccccc}
a_{11} & a_{21} & 0 & 0 & 0
\\
a_{12} & a_{22} & 0 & 0 & 0
\\
0 & 0 & \overline{a}_{11} & \overline{a}_{21} & 0
\\
0 & 0 & \overline{a}_{21} & \overline{a}_{22} & 0
\\
b_1 & b_2 & \overline{b}_1 & \overline{b}_2 & c
\end{array}
\!\right)
\left(\!\!
\begin{array}{c}
\mathcal{L}_1
\\
\mathcal{L}_2
\\
\overline{\mathcal{L}}_1
\\
\overline{\mathcal{L}}_2
\\
\mathcal{T}
\end{array}
\!\!\right),
\]
{\em for some seven local functions:}
\[
\aligned
c\colon\ \ \
M
&
\,\longrightarrow\,
\R\backslash\{0\},
\\
a_{11},\,a_{12},\,a_{21},\,a_{22},\,b_1,\,b_2\,\,\colon\ \ \
M
&
\,\longrightarrow\,
\C,
\endaligned
\]
{\em with:}
\[
0
\,\neq\,
\det\,
\left(\!
\begin{array}{cc}
a_{11} & a_{12} 
\\
a_{21} & a_{22}
\end{array}
\!\right)(q)
\ \ \ \ \ \ \ \ \ \ \ \ \
{\scriptstyle{(\forall\,q\,\in\,M)}}.
\]

\medskip

One easily verifies that the set of matrices:
\[
\!\!\!\!\!\!\!\!\!\!\!\!\!\!\!\!\!\!\!\!
\!\!\!\!\!\!\!\!\!\!\!\!\!\!\!
{\sf G}_{\text{\sf IV}_{\text{\sf 1}}}^{\sf initial}
\,:=\,
\left\{
\left(\!\!
\begin{array}{ccccc}
{\sf a}_{11} & {\sf a}_{21} & 0 & 0 & 0
\\
{\sf a}_{12} & {\sf a}_{22} & 0 & 0 & 0
\\
0 & 0 & \overline{\sf a}_{11} & \overline{\sf a}_{21} & 0
\\
0 & 0 & \overline{\sf a}_{12} & \overline{\sf a}_{22} & 0
\\
{\sf b}_1 & {\sf b}_2 & \overline{\sf b}_1 & \overline{\sf b}_2 & {\sf c}
\end{array}
\!\!\right)
\,\in\,
\mathcal{M}_{5\times 5}(\C)\,\colon\,\,
{\sf c}\,\in\,\R\backslash\{0\},\,\,
{\sf a}_{11},\,{\sf a}_{12},\,{\sf a}_{21},\,{\sf a}_{22},\,
{\sf b}_1,\,{\sf b}_2\,\,
\in\,\C
\right\}
\]
is a (closed) $13$-dimensional
real matrix subgroup of the full:
\[
{\sf GL}_5(\C)
\,:=\,
\left\{
\pi
=
\left(\!\!
\begin{array}{ccccc}
\pi_{1,1} & \pi_{1,2} & \pi_{1,3} & \pi_{1,4} & \pi_{1,5}
\\
\pi_{2,1} & \pi_{2,2} & \pi_{3,3} & \pi_{2,4} & \pi_{2,5}
\\
\pi_{3,1} & \pi_{3,2} & \pi_{3,3} & \pi_{3,4} & \pi_{3,5}
\\
\pi_{4,1} & \pi_{4,2} & \pi_{4,3} & \pi_{4,4} & \pi_{4,5}
\\
\pi_{5,1} & \pi_{5,2} & \pi_{5,3} & \pi_{5,4} & \pi_{5,5}
\end{array}
\!\!\right)
\,\in\,
\mathcal{M}_{5\times 5}(\C)\,\colon\,\,
0
\neq
\det\,\pi
\right\}.
\]

\medskip\noindent{\bf Proposition.}
{\em On a $5$-dimensional hypersurface submanifold:}
\[
\Big(
M^5
\,\subset\,
\C^4
\Big)
\,\,\in\,\,
\text{\sf General Class $\text{\sf IV}_{\text{\sf 1}}$},
\]
{\em having biholomorphically invariant $(1, 0)$ CR bundle:}
\[
T^{1,0}M
\,\subset\,
\C\otimes_\R TM,
\]
{\em for any choice of a pair of local vector field generators:}
\[
\big\{
\mathcal{L}_1,\,\mathcal{L}_2
\big\}
\ \ \ \ \
\text{\rm for}\ \
T^{1,0}M,
\]
{\em the associated frame:}
\[
\big\{
\mathcal{L}_1,\,\mathcal{L}_2,\,
\overline{\mathcal{L}}_1,\,\overline{\mathcal{L}}_2,\,
\isqrt\big[\mathcal{L}_1,\,\overline{\mathcal{L}}_1\big]
\big\}
\,=:\,
\big\{
\mathcal{L}_1,\,\mathcal{L}_2,\,
\overline{\mathcal{L}}_1,\,\overline{\mathcal{L}}_2,\,
\mathcal{T}
\big\}
\]
{\em for the full complexified tangent bundle:}
\[
\C\otimes_\R TM
\]
{\em performs a reduction of the full ${\sf GL}_5 ( \C)$-structure
of $\C \otimes_\R TM$ to the $13$-dimensional subgroup:}
\[
\!\!\!\!\!\!\!\!\!\!\!\!\!\!\!\!\!\!\!\!
\!\!\!\!\!\!\!\!\!\!\!\!\!\!\!\!\!\!\!\!
{\sf G}_{\text{\sf IV}_{\text{\sf 1}}}^{\sf initial}
\,:=\,
\left\{
\left(\!\!
\begin{array}{ccccc}
{\sf a}_{11} & {\sf a}_{21} & 0 & 0 & 0
\\
{\sf a}_{12} & {\sf a}_{22} & 0 & 0 & 0
\\
0 & 0 & \overline{\sf a}_{11} & \overline{\sf a}_{21} & 0
\\
0 & 0 & \overline{\sf a}_{12} & \overline{\sf a}_{22} & 0
\\
{\sf b}_1 & {\sf b}_2 & \overline{\sf b}_1 & \overline{\sf b}_2 & {\sf c}
\end{array}
\!\!\right)
\,\in\,
\mathcal{M}_{5\times 5}(\C)\,\colon\,\,
{\sf c}\,\in\,\R\backslash\{0\},\,\,
{\sf a}_{11},\,{\sf a}_{12},\,{\sf a}_{21},\,{\sf a}_{22},\,
{\sf b}_1,\,{\sf b}_2\,\,
\in\,\C
\right\}.
\qed
\]


\bigskip

\section{\sf General class $\text{\sf IV}_{\text{\sf 2}}$}
\label{general-class-IV-2}
\HEAD{\ref{general-class-IV-2}.~General class 
$\text{\sf IV}_{\text{\sf 2}}$}{
Jo\"el {\sc Merker}, D\'epartement de Math\'ematiques d'Orsay}

\medskip

Equip $\C^3$ with coordinates:
\[
(z_1,z_2,w)
\,\in\,
\C^3.
\]

Let a connected hypersurface:
\[
M^5
\,\subset\,
\C^3
\]
be of smoothness:

\[
\mathcal{C}^\kappa\ \
(\kappa\geqslant 4),
\ \ \ \ \ \ \ \ \ \ \ \ \ \ \ \ \ \ \ \ \
\text{\rm or}
\ \ \ \ \ \ \ \ \ \ \ \ \ \ \ \ \ \ \ \ \
\mathcal{C}^\infty,
\ \ \ \ \ \ \ \ \ \ \ \ \ \ \ \ \ \ \ \ \
\text{\rm or}
\ \ \ \ \ \ \ \ \ \ \ \ \ \ \ \ \ \ \ \ \
\mathcal{C}^\omega.
\]

Pick a point:
\[
p\in M,
\]
and take a (small) open neighborhood:
\[
p
\in
{\sf U}_p
\subset
\C^3.
\]
Take a local frame:
\[
\big\{
\mathcal{L}_1,\,\mathcal{L}_2
\big\}
\ \ \ \ \
\text{\rm for}\ \
T^{1,0}\big(M\cap{\sf U}_p\big).
\]

Recall (\cite{ Merker-Pocchiola-Sabzevari-5-CR-II,
Merker-5-CR-III}) that when:
\[
\Big(
M^5
\,\subset\,
\C^3
\Big)
\,\,\in\,\,
\text{\sf General Class $\text{\sf IV}_{\text{\sf 2}}$},
\]
the kernel, in $T^{1, 0}M$, of the Levi form is an
everywhere of rank $1$ complex subbundle:
\[
K^{1,0}M
\,\subset\,
T^{1,0}M
\]
hence has a vector field generator denoted:
\[
\mathcal{K},
\]
so that, renumbering if necessary, it is more natural
to take:
\[
\big\{\mathcal{K},\,\mathcal{L}_1\big\}
\]
as a local frame for $T^{1, 0} M$.

Then without loss of generality:
\[
\big\{
\mathcal{K},\,\mathcal{L}_1,\,
\overline{\mathcal{K}},\,\overline{\mathcal{L}}_1,\,
\mathcal{T}
\big\}
\]
constitutes a (local) frame for $\C \otimes_\R TM$, 
setting as before:
\[
\mathcal{T}
\,:=\,
\isqrt\,\big[\mathcal{L}_1,
\overline{\mathcal{L}}_1\big].
\]

Because one knows (\cite{ Merker-Pocchiola-Sabzevari-5-CR-II}) that:
\[
\aligned
\big[\mathcal{K},\overline{\mathcal{L}}_1\big]
&
\,\equiv\,
0\ \ \ \ \
\mod\,
\big(
T^{1,0}M
\oplus
T^{0,1}M
\big),
\\
\big[\mathcal{K},\overline{\mathcal{K}}\big]
&
\,\equiv\,
0\ \ \ \ \
\mod\,
\big(
K^{1,0}M
\oplus
K^{0,1}M
\big),
\endaligned
\]
alltogether the brackets are:
\[
\aligned
\isqrt\,\big[\mathcal{L}_1,\overline{\mathcal{L}}_1\big]
&
=
\mathcal{T},
\\
\isqrt\,\big[\mathcal{K},\overline{\mathcal{L}}_1\big]
&
=
A\cdot\mathcal{K}
+
B\cdot\mathcal{L}_1
+
C\cdot\overline{\mathcal{K}}
+
d\cdot\overline{\mathcal{L}}_1,
\\
\isqrt\,\big[\mathcal{L}_1,\overline{\mathcal{K}}\big]
&
=
\overline{C}\cdot\mathcal{K}
+
\overline{D}\cdot\mathcal{L}_1
+
\overline{A}\cdot\overline{\mathcal{K}}
+
\overline{B}\cdot\overline{\mathcal{L}}_1,
\\
\isqrt\,\big[\mathcal{K},\overline{\mathcal{K}}\big]
&
=
E\cdot\mathcal{K}
+
\overline{E}\cdot\overline{\mathcal{K}}
\endaligned
\]
for certain local functions:
\[
A,\,B,\,C,\,D,\,E\,\,\colon\ \ \
M
\,\longrightarrow\,
\C.
\]

\smallskip

Next, take any (local) biholomorphism:
\[
h\colon\ \ \
{\sf U}_p
\overset{\sim}{\,\longrightarrow\,}
{\sf U}_{p'}'
=
h({\sf U}_p)
\ \ \ \ \ \ \ \ \ \ \ \ \
{\scriptstyle{(p'\,=\,\,h(p))}},
\]
which, when ${\sf U}_p$ is small enough, certainly 
transfers $M \cap {\sf U}_p$ to a certain hypersurface:
\[
{M'}^5
\,:=\,
h\big(M\cap{\sf U}_p\big)
\,\subset\,
{\C'}^3,
\]
having the same smoothness (mental exercise):
\[
\mathcal{C}^\kappa\ \
(\kappa\geqslant 4),
\ \ \ \ \ \ \ \ \ \ \ \ \ \ \ \ \ \ \ \ \
\text{\rm or}
\ \ \ \ \ \ \ \ \ \ \ \ \ \ \ \ \ \ \ \ \
\mathcal{C}^\infty,
\ \ \ \ \ \ \ \ \ \ \ \ \ \ \ \ \ \ \ \ \
\text{\rm or}
\ \ \ \ \ \ \ \ \ \ \ \ \ \ \ \ \ \ \ \ \
\mathcal{C}^\omega.
\]
Consider the inverse of $h$:
\[
h'\colon\ \ \
{\sf U}_{p'}'
\overset{\sim}{\,\longrightarrow\,}
{\sf U}_p
\ \ \ \ \ \ \ \ \ \ \ \ \ \ \ \ \ \ \ \ \ \ \ \ \ \ \ \ \ \
{\scriptstyle{(p\,=\,\,h'(p'))}}.
\]

\smallskip

Take similarly a local frame for $TM' \otimes_\R \C$:
\[
\big\{\mathcal{K}',\,\mathcal{L}_1',\,\overline{\mathcal{K}}',\,
\overline{\mathcal{L}}_1',\,\mathcal{T}'
\big\}
\]
with:
\[
\C\mathcal{K}'
=
K^{1,0}M',
\]
and:
\[
\mathcal{T}'
:=
\isqrt\,\big[
\mathcal{L}_1',\overline{\mathcal{L}}_1'
\big].
\]

The invariance of the Levi-kernel:
\[
h_*'\big(K^{1,0}M'\big)
=
K^{1,0}M
\]
yields, as on page~83 of~\cite{ Merker-Pocchiola-Sabzevari-5-CR-II}:
\[
h_*'\big(\mathcal{K}'\big)
=
c\cdot\mathcal{K},
\]
for a certain local function:
\[
c\colon\ \ \ 
M
\,\longrightarrow\,
\C\backslash\{0\},
\]
while of course:
\[
h_*'\big(\mathcal{L}_1'\big)
=
a\cdot\mathcal{L}_1
+
b\cdot\mathcal{K},
\]
for two certain functions:
\[
\aligned
a\colon\ \ \
M
&
\,\longrightarrow\,
\C\backslash\{0\},
\\
b\colon\ \ \
M
&
\,\longrightarrow\,
\C.
\endaligned
\]

\medskip

Now, compute:
\[
\aligned
h_*'\big(\mathcal{T}'\big)
&
=
h_*'\big(\isqrt\,\big[\mathcal{L}_1',\overline{\mathcal{L}}_1'\big]\big)
\\
&
=
\isqrt\,\big[
h_*'\big(\mathcal{L}_1'\big),\,
h_*'\big(\overline{\mathcal{L}}_1'\big)
\big]
\\
&
=
\isqrt\,
\big[
a\,\mathcal{L}_1+b\,\mathcal{K},\,\,
\overline{a}\,\overline{\mathcal{L}}_1
+
\overline{b}\,\overline{\mathcal{K}}
\big]
\\
&
=
a\overline{a}\,\isqrt\,
\big[
\mathcal{L}_1,\overline{\mathcal{L}}_1
\big]
+
b\overline{a}\,\isqrt\,
\big[\mathcal{K},\overline{\mathcal{L}}_1\big]
+
a\overline{b}\,\isqrt\,
\big[\mathcal{L}_1,\overline{\mathcal{K}}\big]
+
b\overline{b}\,\isqrt\,
\big[\mathcal{K},\overline{\mathcal{K}}\big]
+
\\
&
\ \ \ \ \
+
\isqrt\,\big(
a\mathcal{L}_1\big(\overline{a}\big)
+
b\,\mathcal{K}\big(\overline{a}\big)
\big)
\cdot
\overline{\mathcal{L}}_1
+
\isqrt\,\big(
a\mathcal{L}_1\big(\overline{b}\big)
+
b\,\mathcal{K}\big(\overline{b}\big)
\big)
\cdot
\overline{\mathcal{K}}
-
\\
&
\ \ \ \ \
-\,
\isqrt\,\big(
\overline{a}\overline{\mathcal{L}}_1(a)
-
\overline{b}\overline{\mathcal{K}}(a)
\big)
-
\isqrt\,\big(
\overline{a}\overline{\mathcal{L}}_1(b)
-
\overline{b}\overline{\mathcal{K}}(b)
\big),
\endaligned
\]
which is:
\[
\aligned
h_*'\big(\mathcal{T}'\big)
&
=
\big(a\overline{a}\big)
\cdot
\mathcal{T}
+
\\
&
\ \ \ \ \
+
\Big(
b\overline{a}D
+
a\overline{b}\overline{B}
+
\isqrt\,a\mathcal{L}_1\big(\overline{a}\big)
+
\isqrt\,b\mathcal{K}\big(\overline{a}\big)
\Big)
\cdot
\overline{\mathcal{L}}_1
+
\\
&
\ \ \ \ \
+
\Big(
b\overline{a}C
+
a\overline{b}\overline{A}
+
b\overline{b}\overline{E}
+
\isqrt\,a\mathcal{L}_1\big(\overline{b}\big)
+
\isqrt\,b\mathcal{K}\big(\overline{b}\big)
\Big)
\cdot
\overline{\mathcal{K}}
+
\\
&
\ \ \ \ \
+
\Big(
b\overline{a}B
+
a\overline{b}\overline{D}
-
\isqrt\,\overline{a}\overline{\mathcal{L}}_1(a)
-
\isqrt\,\overline{b}\overline{\mathcal{K}}(a)
\Big)
\cdot
\mathcal{L}_1
+
\\
&
\ \ \ \ \
+
\Big(
b\overline{a}A
+
a\overline{b}\overline{C}
+
b\overline{b}E
-
\isqrt\,\overline{a}\overline{\mathcal{L}}_1(b)
-
\isqrt\,\overline{b}\overline{\mathcal{K}}(b)
\Big)
\cdot
\mathcal{K}.
\endaligned
\]

Two coefficient-function appear (together with their
conjugates):
\[
\aligned
d
&
:=
b\overline{a}B
+
a\overline{b}\overline{D}
-
\isqrt\,\overline{a}\overline{\mathcal{L}}_1(a)
-
\isqrt\,\overline{b}\overline{\mathcal{K}}(a),
\\
e
&
:=
b\overline{a}A
+
a\overline{b}\overline{C}
+
b\overline{b}E
-
\isqrt\,\overline{a}\overline{\mathcal{L}}_1(b)
-
\isqrt\,\overline{b}\overline{\mathcal{K}}(b),
\endaligned
\]
of course complex-valued:
\[
d,\,\,e\,\,\colon\ \ \
M
\,\longrightarrow\,
\C.
\]

\medskip\noindent{\bf Summary.}
{\em Through any local biholomorphic equivalences
between CR-generic submanifolds of $\C^3$ belonging
to the General Class} $\text{\sf IV}_{\text{\sf 2}}$:
\[
M^5
\overset{\sim}{\,\longrightarrow\,}
{M'}^5,
\]
{\em for any choices of local vector
field generators of the Levi-kernel bundles:}
\[
\aligned
&
\big\{\mathcal{K}\big\}
\ \ \ \ \
\text{\rm for}\ \
T^{1,0}M,
\\
&
\big\{\mathcal{K}'\big\}
\ \ \ \ \
\text{\rm for}\ \
T^{1,0}M',
\endaligned
\]
and for any choice of completion frame:
\[
\aligned
&
\big\{\mathcal{K},\,\mathcal{L}_1\big\}
\ \ \ \ \
\text{\rm for}\ \
T^{1,0}M,
\\
&
\big\{\mathcal{K}',\mathcal{L}_1'\big\}
\ \ \ \ \
\text{\rm for}\ \
T^{1,0}M',
\endaligned
\]
{\em the transfer of frame obeys the rule:}
\[
\left(\!\!
\begin{array}{c}
\mathcal{K}'
\\
\mathcal{L}_1'
\\
\overline{\mathcal{K}}'
\\
\overline{\mathcal{L}}_1'
\\
\mathcal{T}'
\end{array}
\!\!\right)
=
\left(\!
\begin{array}{ccccc}
c & 0 & 0 & 0 & 0
\\
b & a & 0 & 0 & 0
\\
0 & 0 & \overline{c} & 0 & 0
\\
0 & 0 & \overline{b} & \overline{a} & 0
\\
e & d & \overline{e} & \overline{d} & a\overline{a}
\end{array}
\!\right)
\left(\!\!
\begin{array}{c}
\mathcal{K}
\\
\mathcal{L}_1
\\
\overline{\mathcal{K}}
\\
\overline{\mathcal{L}}_1
\\
\mathcal{T}
\end{array}
\!\!\right),
\]
{\em for some five local functions:}
\[
\aligned
a,\,c\colon\ \ \
M
&
\,\longrightarrow\,
\C\backslash\{0\},
\\
b,\,d,\,e\,\colon\ \ \
M
&
\,\longrightarrow\,
\C.\qed
\endaligned
\]

\medskip

One easily verifies that the set of matrices:
\[
\!\!\!\!\!\!\!\!\!\!\!\!\!\!\!\!\!\!\!\!
\!\!\!\!\!\!\!\!\!\!\!\!\!\!\!
{\sf G}_{\text{\sf IV}_{\text{\sf 2}}}^{\sf initial}
\,:=\,
\left\{
\left(\!\!
\begin{array}{ccccc}
{\sf c} & 0 & 0 & 0 & 0
\\
{\sf b} & {\sf a} & 0 & 0 & 0
\\
0 & 0 & \overline{\sf c} & 0 & 0
\\
0 & 0 & \overline{\sf b} & \overline{\sf a} & 0
\\
{\sf e} & {\sf d} & \overline{\sf e} & \overline{\sf d} & 
{\sf a}\overline{\sf a}
\end{array}
\!\!\right)
\,\in\,
\mathcal{M}_{5\times 5}(\C)\,\colon\,\,
{\sf a},\,{\sf c}\,\in\,\C\backslash\{0\},\,\,
{\sf b},\,{\sf d},\,{\sf e}
\,\in\,\C
\right\}
\]
is a (closed) $10$-dimensional
real matrix subgroup of the full:
\[
{\sf GL}_5(\C)
\,:=\,
\left\{
\pi
=
\left(\!\!
\begin{array}{ccccc}
\pi_{1,1} & \pi_{1,2} & \pi_{1,3} & \pi_{1,4} & \pi_{1,5}
\\
\pi_{2,1} & \pi_{2,2} & \pi_{3,3} & \pi_{2,4} & \pi_{2,5}
\\
\pi_{3,1} & \pi_{3,2} & \pi_{3,3} & \pi_{3,4} & \pi_{3,5}
\\
\pi_{4,1} & \pi_{4,2} & \pi_{4,3} & \pi_{4,4} & \pi_{4,5}
\\
\pi_{5,1} & \pi_{5,2} & \pi_{5,3} & \pi_{5,4} & \pi_{5,5}
\end{array}
\!\!\right)
\,\in\,
\mathcal{M}_{5\times 5}(\C)\,\colon\,\,
0
\neq
\det\,\pi
\right\}.
\]

\medskip\noindent{\bf Proposition.}
{\em On a $5$-dimensional hypersurface submanifold:}
\[
\Big(
M^5
\,\subset\,
\C^4
\Big)
\,\,\in\,\,
\text{\sf General Class $\text{\sf IV}_{\text{\sf 2}}$},
\]
{\em having biholomorphically invariant $(1, 0)$ CR bundle:}
\[
T^{1,0}M
\,\subset\,
\C\otimes_\R TM,
\]
{\em for any choice of a pair of vector field generators:}
\[
\aligned
&
\big\{
\mathcal{K}
\big\}
\ \ \ \ \ \ \ \ \ \ \ \ \,
\text{\rm for}\ \
K^{1,0}M,
\\
&
\big\{
\mathcal{K},\,\mathcal{L}_1
\big\}
\ \ \ \ \
\text{\rm for}\ \
T^{1,0}M,
\endaligned
\]
{\em the associated frame:}
\[
\big\{
\mathcal{K},\,\mathcal{L}_1,\,
\overline{\mathcal{K}},\,\overline{\mathcal{L}}_1,\,
\isqrt\big[\mathcal{L}_1,\,\overline{\mathcal{L}}_1\big]
\big\}
\,=:\,
\big\{
\mathcal{K},\,\mathcal{L}_1,\,
\overline{\mathcal{K}},\,\overline{\mathcal{L}}_1,\,
\mathcal{T}
\big\}
\]
{\em for the full complexified tangent bundle:}
\[
\C\otimes_\R TM
\]
{\em performs a reduction of the full ${\sf GL}_5 ( \C)$-structure
of $\C \otimes_\R TM$ to the $10$-dimensional subgroup:}
\[
\!\!\!\!\!\!\!\!\!\!\!\!\!\!\!\!\!\!\!\!
\!\!\!\!\!\!\!\!\!\!\!\!\!\!\!\!\!\!\!\!
{\sf G}_{\text{\sf IV}_{\text{\sf 2}}}^{\sf initial}
\,:=\,
\left\{
\left(\!\!
\begin{array}{ccccc}
{\sf c} & 0 & 0 & 0 & 0
\\
{\sf b} & {\sf a} & 0 & 0 & 0
\\
0 & 0 & \overline{\sf c} & 0 & 0
\\
0 & 0 & \overline{\sf b} & \overline{\sf a} & 0
\\
{\sf e} & {\sf d} & \overline{\sf e} & \overline{\sf d} & 
{\sf a}\overline{\sf a}
\end{array}
\!\!\right)
\,\in\,
\mathcal{M}_{5\times 5}(\C)\,\colon\,\,
{\sf a},\,{\sf c}\,\in\,\C\backslash\{0\},\,\,
{\sf b},\,{\sf d},\,{\sf e}
\,\in\,\C
\right\}.
\qed
\]


\vfill\end{document}

%% file: macros.tex
\newcommand{\C}{\mathbb{C}}

\newcommand{\R}{\mathbb{R}}

\newcommand{\red}{\textcolor{red}}

\newcommand{\HEAD}[2]{%
\pagestyle{fancy}
\fancyhead[RO]{\tiny\sf\thepage}
\fancyhead[CO]{{\tiny\sf #1}}
\fancyhead[LE]{\tiny\sf\thepage}
\fancyhead[CE]{{\tiny\sf #2}}
\fancyfoot{}}

\theoremstyle{definition}

\renewcommand{\det}{\text{\footnotesize\sf det}}

\newcommand{\isqrt}{{\scriptstyle{\sqrt{-1}}}}

\renewcommand{\lim}{\text{\footnotesize\sf lim}}

\renewcommand{\mod}{\text{\footnotesize\sf mod}}

\newcommand{\rank}{\text{\footnotesize\sf rank}}

\newcommand{\vf}{\vfill\end{document}}

\newcommand{\zero}[1]{\underline{#1}_{\red{\circ}}}
